\documentclass[11pt]{article}
\usepackage{amssymb}
\usepackage{amsmath}
\usepackage{amsthm}
\usepackage{latexsym}
\usepackage{graphicx}
\usepackage{enumerate}

\numberwithin{equation}{section}

\def\rc#1{\frac{1}{#1}}
\def\qed{\square}

\def\pmatrix#1{\left(\begin{matrix} #1 \end{matrix}\right)}

\newtheorem{Lemma}{Lemma}
\newtheorem{Corollary}{Corollary}

\newtheorem{Theorem}{Theorem}
\theoremstyle{definition}
\newtheorem*{rmk}{Remark}

\newtheorem*{defn}{Definition}
\newcommand{\ba}[1]{\begin{array}{@{}#1@{}}}
\newcommand{\ea}{\end{array}}
\long\def\symbolfootnote[#1]#2{\begingroup%
\def\thefootnote{\fnsymbol{footnote}}\footnote[#1]{#2}\endgroup}

\def\emx{{\cal E}}
\def\jmx{{\cal J}}
\def\imx{{\cal I}}
\def\omx{{\cal O}}
\def\xmx{{\cal X}}

\begin{document}
\newcommand{\ul}{\underline}
\newcommand{\be}{\begin{equation}}
\newcommand{\ee}{\end{equation}}
\newcommand{\ben}{\begin{enumerate}}
\newcommand{\een}{\end{enumerate}}

\long\def\symbolfootnote[#1]#2{\begingroup%
\def\thefootnote{\fnsymbol{footnote}}\footnote[#1]{#2}\endgroup}
\baselineskip = 15pt
\noindent

\title{On Superalgebras of Matrices with Symmetry Properties}
\date{\today}

\author{S. L. Hill, M. C. Lettington\\ and K. M. Schmidt (Cardiff)$$}
\maketitle
\begin{abstract}
It is known that semi-magic square matrices form a 2-graded algebra
or superalgebra with the even and odd subspaces under centre-point reflection
symmetry as the two components. We show that other symmetries which have been
studied for square matrices give rise to similar superalgebra structures,
pointing to novel symmetry types in their complementary parts. In particular,
this provides a unifying framework for the composite `most perfect square'
symmetry and the related class of `reversible squares';
moreover, the semi-magic square algebra is identified as part of a 2-gradation
of the general square matrix algebra.
We derive explicit representation formulae for matrices of all symmetry types
considered, which can be used to construct all such matrices.
\end{abstract}

\section{Matrix Symmetry Type Spaces}
\label{sMSTS}

In this paper we present a novel approach to the classification of particular families of $n\times n$ matrices, defined by their symmetry properties, in terms of $\mathbb{Z}_2$ graded algebras.
The latter type of algebra (also known as \emph{superalgebra}) has a decomposition into an `even' subalgebra and an `odd' complementary part which is a bimodule over the `even' subalgebra and squares into it.
The families of matrices considered here are derived from nine fundamental symmetry properties that generate corresponding matrix symmetry vector spaces.
Using a block matrix representation introduced in \cite{hill}, we find that these matrix spaces arrange into four $\mathbb{Z}_2$ graded algebras and a single algebra (the type R algebra defined below).

These algebraic structures enable us to analyse more specialised algebras of matrices, defined by compositions of
these symmetry properties. Such matrix families encompass some well known symmetry types such as the sets of semimagic square matrices \cite{jesus}, the associated magic square matrices \cite{andrewsmagic}, most perfect square matrices, and the reversible square matrices of $\cite{bree}$. Accompanying these matrix families in their respective $\mathbb{Z}_2$ graded algebras, we find hitherto undocumented matrix symmetry types such as types N, Q and V defined in this paper.
The present findings build on recent work \cite{brunnock, hill, mcl1, mcl2} to provide insight into the deeper algebraic structures underpinning an area of mathematics that has been of interest for many years.
In the process, we derive matrix algebraic characterisations of the symmetries and find representation formulae for the matrices of each type, which can be used to construct them.

We begin with a review of known symmetry types, from which our basic symmetries are then distilled.

Semimagic squares are defined by the property that all rows and all columns add
up to the same constant. In an associated semimagic square, opposite entries
with respect to the centre of the square also add up to the same number.
A balanced semimagic square has the complementary property that the opposite
entries are equal, so the square has a half-turn rotational symmetry.

Most perfect squares are semimagic squares with the additional properties that
all $2\times 2$ blocks of numbers add up to the same constant and that the
square has the strong pandiagonal property, so all pairs of entries half the
size of the square apart along a general diagonal (i.e. any line parallel to either of the two main diagonals) add up to the same constant. Clearly
this definition only makes sense for squares of even dimension.

Reversible squares are defined by the properties that all pairs of entries
on a row or column which have the same distance from the centre add up to the
same constant, and that for any rectangular submatrix, the two pairs of
diagonally opposite vertex entries add up to the same constant.

It was shown in \cite{brunnock} that, after the removal of a common
two-sided ideal, the semimagic squares, considered as square matrices with the
usual matrix addition and multiplication operations, form an algebra which
has the form of a $\mathbb{Z}_2$-graded algebra, with the balanced semimagic
matrices as `even' subalgebra and the associated magic matrices as `odd'
complementary direct summand.
It is our aim in the present paper to explore the algebraic properties of
the other types of matrices mentioned above, and of more general symmetry types
arising in their definitions.
Along the way, we also establish representation formulae for matrices of different
symmetry types, which make their algebraic behaviour more transparent and also
provide a simple way of constructing matrices of a particular type.
It turns out that the $\mathbb{Z}_2$-graded
algebra structure recurs in various guises.

The set of most perfect square matrices and the set of reversible
square matrices do not themselves form subalgebras of the general algebra of
square matrices, similarly to the case of the associated semimagic square
matrices. However, we identify suitable complementing
subalgebras, analogous to the set of balanced semimagic square matrices, which
extend these sets to $\mathbb{Z}_2$-graded algebras.
For this purpose, we need to break down the composite definitions into a
number of simpler symmetry conditions.

We shall consider the following symmetries of a square matrix
$M=\left ( M_{i,j}\right )_{i,j \in \mathbb{Z}_n} \in \mathbb{R}^{n\times n}$.
Note that the indices are considered to be elements of the cyclic ring
$\mathbb{Z}_n := \mathbb{Z}/{n \mathbb{Z}}$,
and all calculations with indices are performed in this ring, i.e. modulo $n$.
The top left corner of the matrix will have indices
$(1,1) \in \mathbb{Z}_n^2$.

\begin{description}
\item{(S)} Semimagic property of weight $w$:
\[
\sum_{j\in\mathbb{Z}_n} M_{i,j}=\sum_{j\in\mathbb{Z}_n} M_{j,i}=nw
 \qquad (i \in \mathbb{Z}_n).
\]
\item{(A)} Associated property of weight $w$ : $M_{i,j} + M_{n+1-i, n+1-j} = 2 w$
$(i, j \in \mathbb{Z}_n)$.
\item{(B)} Balanced property: $M_{i,j} - M_{n+1-i, n+1-j} = 0$
$(i, j \in \mathbb{Z}_n)$.
%\item{(D)} Diagonal sum symmetry of weight $w$:
%\[
%\sum_{i\in\mathbb{Z}_n} m_{i,i}=\sum_{i\in\mathbb{Z}_n} m_{i, n+1-i}=nw.
%\]
%\item{(W)} Weak pandiagonal symmetry of weight $w$:
%\[
%\sum_{j \in\mathbb{Z}_n} m_{q+j,j}=\sum_{j \in\mathbb{Z}_n} m_{q-j,j}=nw
%\qquad (q \in \mathbb{Z}_n).
%\]
\item{(R)} Row and column reverse property:
\begin{align*}
 M_{i,j} +M_{i,n+1-j}&=M_{i,k} +M_{i,n+1-k}, \\
 M_{i,j}+M_{n+1-i,j}&=M_{k,j}+M_{n+1-k,j}
 \qquad (i,j,k \in \mathbb{Z}_n).
\end{align*}
\item{(V)} Vertex cross sum property: $M_{i,j}+M_{k,l}=M_{i,l}+M_{k,j}$
$(i,j,k,l \in \mathbb{Z}_n)$.
\end{description}
\noindent
In the case where $n = 2 \nu$ is even, we also consider the following symmetries.
\begin{description}
\item{(M)} $2\times 2$ array sum property of weight $w$: $M_{i,j} +M_{i,j+1} + M_{i+1,j}+M_{i+1,j+1} =4w$
$(i,j \in \mathbb{Z}_n)$
\hfill\break
{\bf and}
alternating sum property:
\[
 \sum_{i,j \in \mathbb{Z}_n} (-1)^{i+j} M_{i,j} = 0.
\]
\item{(N)} Consecutive row and column alternating sum property:
\[
\sum_{i\in\mathbb{Z}_n} (-1)^i (M_{i,j} + M_{i,{j+1}}) =
\sum_{i\in\mathbb{Z}_n} (-1)^i (M_{j,i} + M_{{j+1},i}) = 0
\qquad (j \in \mathbb{Z}_n).
\]
\item{(P)} Strong pandiagonal property of weight $w$: $M_{i,j} + M_{i+\nu,j+\nu}=2w$
$(i,j \in \mathbb{Z}_n)$.
\item{(Q)} Quartered sum property: $M_{i,j} - M_{i+\nu,j+\nu}=0$
$(i,j \in \mathbb{Z}_n)$.
\end{description}

\noindent
{\it Remarks.}
1. The following two further symmetries are often considered. The first is the
property that both of the two diagonals of a semimagic square matrix of weight $w$ each
add up to $n w$; this is then called a {\it magic square matrix}.
This property evidently follows from (A) or (P).
The second is the (weak) pandiagonal property, where all general (cyclically broken)
diagonals of the matrix add up to $n w$; this clearly follows from (P).
We do not consider these two symmetries in this paper, except where they
naturally follow from stronger properties.

2. Property (M) does not at face value presuppose even matrix dimension $n$.
However, if $n$ is odd, then only the null matrix
\begin{equation}
\label{eomxdef}
\omx_n = (0)_{i,j=1}^n
\end{equation}
has this property; see Lemma \ref{lnoModd} below.
We note that there are odd-dimensional non-trivial matrices with property (N), e.g.
\begin{equation*}
 M = \pmatrix{1 & 2 & 1 \cr 1 & 0 & -1 \cr 0 & -2 & -2}.
\end{equation*}

3. A most perfect square matrix has properties (M), (P) and (S).
Note that in the original definition by Ollerenshaw (see \cite{bree} page 12), the alternating sum
property part of (M) was not stipulated, however it is
already implied by (P) in the case of even dimension $n = 2 \nu$;
indeed, then
\begin{align*}
\sum_{i,j\in\mathbb{Z}_n} (-1)^{i+j} M_{i,j}
&= \sum_{i,j=1}^\nu \left((-1)^{i+j} M_{i,j} + (-1)^{i+j+2\nu} M_{i+\nu,j+\nu}
   + (-1)^{i+j+\nu} (M_{i+\nu,j} + M_{i,j+\nu})\right) \\
&= \sum_{i,j=1}^\nu (-1)^{i+j} (2 w + (-1)^\nu 2 w) = 0
\end{align*}
both for even and odd $\nu$.
When property (M) is considered by itself, the additional alternating sum property
is essential to give a clear separation from property (N), see Theorem \ref{tdisjoint}.

4. Reversible squares have properties (R) and (V). Moreover, as we shall see below
in Corollary \ref{cdirsum} and Theorem \ref{tSValg},
property (V) somewhat surprisingly also plays a role as a complement to property
(S). Reversible squares arose from Ollerenshaw's adaptation of a 1938 construction \cite{diabolic} of
Rosser and
Walker, which she used to enumerate the number of doubly-even order most perfect square matrices.

5. Properties (N) and (Q) have not previously been studied; we identify them here
as natural complements to properties (M) and (P), respectively.
In the case of (Q) this is easy to understand;
a matrix with symmetry (Q) is of the form
\begin{equation}
\label{eQstruct}
\left(\begin{matrix} A & B \\ B & A \end{matrix}\right)
\end{equation}
with $A, B \in\mathbb{R}^{\nu\times\nu}$, resembling a quartered shield in heraldry,
whereas matrices with symmetry (P) have the structure
\begin{equation}
\label{ePstruct}
\left(\begin{matrix} A & B \\ -B & -A \end{matrix}\right)
\end{equation}
with $A, B \in \mathbb{R}^{\nu\times\nu}$,
so by a straightforward calculation any $n\times n$ matrix can be written as a sum of
type (P) and type (Q) matrices.
Symmetry (N), which means that the alternating sum of each row is the negative of the
alternating sum of a neighbouring row, and similarly for the columns,
is not very intuitive; it arises
as a complementary property to (M) shared by products of matrices having property (M), see Corollary \ref{cdirsum} and Theorem \ref{tNMalg}.

\begin{Lemma}
\label{lnoModd}
Let $n\in\mathbb{N}$ be odd and $M \in \mathbb{R}^{n\times n}$ a matrix with
property (M). Then $M = \omx_n$.
\end{Lemma}

\medskip\noindent
{\it Proof.\/}
By the $2\times 2$ array sum property,
we have for each $i \in \mathbb{Z}_n$
\begin{align*}
 M_{i,1} + M_{i+1,1} &= 4w - M_{i,2} - M_{i+1,2}
= M_{i,3} + M_{i+1,3} = 4w - M_{i,4} - M_{i+1,4} \\
&= \cdots = M_{i,n} + M_{i+1,n} = 4w - M_{i,1} - M_{i+1,1},
\end{align*}
so $M_{i,j} + M_{i+1,j} = 2 w$ for all $i, j \in \mathbb{Z}_n$.
Hence, for each $j \in\mathbb{Z}_n$,
\begin{equation*}
 M_{1,j} = 2w - M_{2,j} = M_{3,j} = \dots = M_{n,j} = 2w - M_{1,j},
\end{equation*}
which implies $M_{i,j} = w$ $(i,j\in\mathbb{Z}_n)$.
But then the alternating sum property requires $w = 0$.
\phantom{.}\hfill$\qed$\par\noindent

\bigskip\noindent
In \cite{brunnock, hill} it was observed that the matrix
\begin{equation}
\label{eemxdef}
\emx_n = (1)_{i,j=1}^n
\end{equation}
generates a two-sided ideal in the algebra
of matrices having property (S), and that this ideal is the intersection of the
subspace of matrices with properties (A), (S) and the subalgebra of matrices
with properties (B), (S). Also, subtracting $w\emx_n$ from the matrices which
have (A) with weight $w$ gives weight $0$ matrices with the same symmetry.

In fact, the very simple matrix $\emx_n$ shares all of the above symmetries.

\begin{Lemma}
\label{lemxprop}
Let $n \in\mathbb{N}$. The matrix $\emx_n$ has properties (S), (A), (B), (R), (V),
and, if $n$ is even, also (M), (N), (P) and (Q); where applicable, its weight is
$w = 1$.
\end{Lemma}
\noindent
In consequence, we can often restrict our attention to the weightless case
$w = 0$ by subtracting a suitable multiple of $\emx_n$ from the matrices under
consideration; we shall do this regularly with properties (A), (M) and (P).

Furthermore, all of the above symmetry properties are linear (with weight
either fixed to $0$ or left variable) and hence give rise to vector spaces
of matrices as follows.
The additional requirement in the definition of $V_n$ corresponds to the restriction
to weight $0$ in $A_n$, $M_n$ and $P_n$.
\begin{defn}
Let $n \in \mathbb{N}$.
We define the following matrix symmetry type spaces.
\begin{description}
\item{} $S_n = \{M \in\mathbb{R}^{n\times n} \mid M \textrm{ has property (S) with some weight $w$}\},$
\item{} $A_n = \{M \in\mathbb{R}^{n\times n} \mid M \textrm{ has property (A) with weight } 0\},$
\item{} $B_n = \{M \in\mathbb{R}^{n\times n} \mid M \textrm{ has property (B)}\},$
\item{} $R_n = \{M \in\mathbb{R}^{n\times n} \mid M \textrm{ has property (R)}\},$
\item{} $V_n = \{M = (M_{i,j})_{i,j=1}^n \in\mathbb{R}^{n\times n} \mid M \textrm{ has property (V), and } \sum_{i,j=1}^n M_{i,j} = 0\}.$
\end{description}
For even $n$, we also define the symmetry type spaces
\begin{description}
\item{} $M_n = \{M \in\mathbb{R}^{n\times n} \mid M \textrm{ has property (M) with weight } 0\},$
\item{} $N_n = \{M \in\mathbb{R}^{n\times n} \mid M \textrm{ has property (N)}\},$
\item{} $P_n = \{M \in\mathbb{R}^{n\times n} \mid M \textrm{ has property (P) with weight } 0\},$
\item{} $Q_n = \{M \in\mathbb{R}^{n\times n} \mid M \textrm{ has property (Q)}\}.$
\end{description}
We also consider intersections of these spaces in Sections \ref{sCSMPS} and \ref{sCSRS}.
\end{defn}
\medskip\noindent
{\it Remark\/} on the history of these symmetry terms. The magic square of order $3$ or `Loh Shu square' was known in China as early as the Warring States period (481 -- 221 BC) \cite{andrewsmagic}. An associated $4\times 4$ magic square is depicted by D\"urer in his 1514 work `Melencolia I', and in the 18th century Euler derived a method of constructing semimagic squares whilst attempting to solve the `six officers problem', which was shown to have no solution by Tarry in 1890 \cite{kraitchik} p.~159.
%(see \symbolfootnote[1]{The six officers problem: 36 officers are chosen, six from each of six different regiments and six from each of six different ranks. They are to be placed in six groups of six officers in such a way that in each group there is one officer from every regiment and one officers of every rank. Euler supsected that there could be no solution to this problem, but he did not succeed in probing this. In 1890 a proof was found by Tarry (see Kraitchik, \cite{kraitchik}~p~159). Euler's search led him to develop his method for constructing semi-magic squares.}
Pandiagonal magic squares were considered in an 1878 paper by Frost, and in 1897, McClintock defined and gave a construction for most perfect magic squares. Hence, in these mixed forms, the type S, A, M and P symmetries were all well established by 1900.
Ollerenshaw's detailed study of the most perfect square (composite (M), (S), (P)) symmetry gave rise to the definition and study \cite{bree} of the reversible square (composite (R), (V)) symmetry. The link between the type A symmetry to the centro-symmetric type B symmetry has been established in \cite{mcl1}.
As far as the authors are aware, the three symmetries (V) (as a separate symmetry), (N) and (Q) are first introduced here to
complement the semimagic, most perfect and pandiagonal symmetries, respectively.

\medskip\noindent
The different symmetry properties were defined above by reference to the individual
matrix entries. This is descriptive and helps visualise each particular matrix
symmetry, but it is rather inconvenient for studying the algebraic properties of
the symmetry type. For this purpose, we now give an equivalent characterisation
of the symmetries in terms of matrix algebra.

Here and in the following we shall use the following notation.
We write
$0_n$ for the null vector in $\mathbb{R}^n$, and $1_n$ for the vector in this space
which has all entries equal to 1.
Moreover, we define the alternating vector $\S_n$ which has $(-1)^{j-1}$ for its
$j$th entry, $j \in \{1, \dots, n\}$.
These and other vectors in $\mathbb{R}^n$ are
considered as column vectors; we denote row vectors by the transpose of column
vectors, writing $v^T$ for a row vector, where $v \in\mathbb{R}^n$.
Thus for even $n\in\mathbb{N}$, we have
\[
 \S_n = (1, -1, 1, -1, \dots, 1, -1)^T \in \mathbb{R}^n,
\]
and this vector will be orthogonal on $1_n$, but this will not be the case if
$n$ is odd, since then $\S_n$ will have 1 as its last entry.

In addition to the matrices
$\emx_n$ and $\omx_n$ already defined in (\ref{eemxdef}) and (\ref{eomxdef}), respectively,
we use the special matrices
$\jmx_n = (\delta_{i, n+1-j})_{i,j = 1}^n \in {\mathbb R}^{n \times n},$
which has entries
$1$
on the antidiagonal and
$0$
otherwise,
and the $n \times n$ unit matrix
$\imx_n = (\delta_{i,j})_{i,j = 1}^n$,
where $\delta_{i,j}$ is the Kronecker symbol.
As usual, we denote by $X^\bot = \{u \in \mathbb{R}^n \mid u^T v = 0 \ (v \in X)\}$
the orthogonal complement of a set $X \subset \mathbb{R}^n$.
\begin{Theorem}
\label{tsymeq}
Let $M \in \mathbb{R}^{n\times n}$, $n \in \mathbb{N}$.
Then
\begin{description}
\item{(a)}
$M \in S_n$ if and only if $1_n^T M u = 0 = u^T M 1_n$ $(u \in \{1_n\}^\bot)$;
\item{(b)}
$M \in A_n$ if and only if $M + \jmx_n M \jmx_n = \omx_n$;
\item{(c)}
$M \in B_n$ if and only if $M = \jmx_n M \jmx_n$;
\item{(d)}
$M \in R_n$ if and only if
$(M + M\jmx_n) u = 0$ and $(M^T + M^T\jmx_n) u = 0$ $(u \in \{1_n\}^\bot)$;
these are equivalent to
$(M^T + \jmx_n M^T)\mathbb{R}^n \subset \mathbb{R} 1_n$ and
$(M + \jmx_n M)\mathbb{R}^n \subset \mathbb{R} 1_n$,
respectively;
\item{(e)}
$M \in V_n$ if and only if
$u^T M v = 0$ $(u, v \in \{1_n\}^\bot)$ and $1_n^T M 1_n = 0$;
\item{(f)}
if $n$ is even, then
$M \in M_n$ if and only if
$u^T M v = 0$ $(u, v \in \{\S_n\}^\bot)$ and $\S_n^T M \S_n = 0$;
\item{(g)}
if $n$ is even, then
$M \in N_n$ if and only if
$\S_n^T M u = 0 = u^T M \S_n$ $(u \in \{\S_n\}^\bot)$.
\end{description}
\end{Theorem}
\noindent
{\it Remark.}
The symmetry types (Q) and (weightless) (P) are conveniently described by their
block matrix structures (\ref{eQstruct}) and (\ref{ePstruct}), respectively.
Note that $S_n$, $V_n$ closely parallel $N_n$ and $M_n$, respectively, with the
vector $\S_n$ taking the role of $1_n$ for the latter pair.

\medskip\noindent
For later use, we also define symmetry spaces $M_n$ and $N_n$ for {\it odd\/}
$n$ in terms of the properties in Theorem \ref{tsymeq} (f), (g); this will
be useful in Theorems \ref{tNblock} and \ref{tMblock}. Note, however, that
the elements of these spaces will {\it not\/} in general have the symmetries
(N) or (M), respectively; in particular, $M_n$ will contain non-trivial matrices
notwithstanding Lemma \ref{lnoModd}.

\medskip\noindent
{\bf Definition.}
Let $n\in\mathbb{N}$ be odd. Then we define the symmetry type spaces
\begin{description}
\item{} $M_n = \{M\in\mathbb{R}^{n\times n}\mid u^T M v = 0 \ (u, v\in\{\S_n\}^\bot),
 \S_n^T M \S_n = 0\}$,
\item{} $N_n = \{M\in\mathbb{R}^{n\times n}\mid \S_n^T M u = 0 = u^T M \S_n
 \ (u\in\{\S_n\}^\bot)\}$.
\end{description}

\bigskip
\noindent
In the proof of Theorem \ref{tsymeq} (a) and later on we shall use the following
observation that the conditions in (a) and (g) are equivalent to an eigenvalue
property of $M$ and $M^T$.

\medskip\noindent
\begin{Lemma}
\label{lMNSVeigen}
Let $n \in \mathbb{N}$ and $y \in \mathbb{R}^n\setminus\{0_n\}$.
Then $M \in \mathbb{R}^{n \times n}$
satisfies
\begin{equation*}
y^T M u = 0 = u^T M y \qquad (u \in \{y\}^\bot)
\end{equation*}
if and only if there is some $\lambda\in\mathbb{R}$
such that $M y = \lambda y$, $M^T y = \lambda y$.
\end{Lemma}

\noindent
{\it Proof.}
Since $0 = u^T M y$ for all $u \in \{y\}^\bot$, we see that
$M y \in \{y\}^{\bot\bot} = \mathbb{R} y$, so there is some $\lambda\in\mathbb{R}$
such that $M y = \lambda y$.
Similarly, $0 = y^T M u = (u^T M^T y)^T$ shows that there is some
$\lambda'\in\mathbb{R}$ such that $M^T y = \lambda' y$.
Hence
\begin{equation*}
\lambda y^T y = y^T M y = (M^T y)^T y = \lambda' y^T y,
\end{equation*}
and as $y^T y \neq 0$, it follows that $\lambda' = \lambda$.
The converse statement is obvious.
\phantom{.}\hfill$\qed$\par\noindent

\medskip\noindent
{\it Proof\/} of Theorem \ref{tsymeq}.
(a) Property (S) can be rewritten in the form $M 1_n = M^T 1_n = n w 1_n$, so
the equivalence follows by Lemma \ref{lMNSVeigen} with $y = 1_n$.
(b) and (c) are straightforward, noting that conjugation with $\jmx_n$ rotates
the matrix by a half-turn.

For (d), note that $(M + M\jmx_n) u = 0$ $(u \in \{1_n\}^\bot)$ means that
$M + M\jmx_n = (M_{i,j} + M_{i,n+1-j})_{i,j\in\mathbb{Z}_n}$ has constant rows.
Also, $(M + \jmx_n M)\mathbb{R}^n \subset \mathbb{R} 1_n$ means that
$M + \jmx_n M = (M_{i,j} + M{n+1-i,j})_{i,j\in\mathbb{Z}_n}$ has constant columns.
These statements are equivalent to (R).
The other equivalent equations follow by considering $M^T$.

For (e), first note that $1_n^T M 1_n = \sum_{i,j=1}^n M_{i,j}$.
Now consider the vectors
$v_j$,
$j \in \{1, n-1\}$,
defined such that
$v_j$
has
$1$
in the $j$th and
$-1$
in the $(j+1)$st positions, and zeros otherwise. These vectors form a basis of
$\{1_{n}\}^\bot;$
indeed, any vector
$u = (u_1, u_2, \dots, u_{n})^T$
such that
$\sum_{k=1}^{n} u_k = 0$
can be rewritten as
\begin{align}
 u &= u_1 v_1 + (u_1 + u_2) v_2 + (u_1 + u_2 + u_3) v_3 + \cdots + (u_1 + u_2 + \cdots
      + u_{n-1}) v_{n-1}.
\nonumber\end{align}
Now, for any
$j, k \in \{1, \dots, n-1\}$,
\begin{align}
 v_j^T M v_k &= M_{j,k} + M_{j+1,k+1} - M_{j,k+1} - M_{j+1,k} = 0
\nonumber\end{align}
by property (V), and (e) follows by bilinearity.
%\par
Conversely, if (e) holds and
$j, k, l, m \in \{1, \dots, n\},$
let
$u$
be the vector such that
$u_j = 1$,
$u_l = -1$
and all other entries vanish, and let
$v$
be the vector such that
$v_k = 1$,
$v_m = -1$
and all other entries vanish.
Then
$u, v \in \{1_{n}\}^\bot$,
so
\begin{align}
 0 = u^T M v &= M_{j,k} + M_{l,m} - M_{j,m} - M_{l,k},
\nonumber\end{align}
and hence
$M$
has property (V).

For (f), note first that $\S_n^T M \S_n = \sum_{i,j=1}^n (-1)^{i+j} M_{i,j}$.
Further, the $2 \times 2$ array sum property with weight 0 can be expressed as
\begin{align}
 v_i^T M v_j &= 0
 \qquad (i, j \in \{1, \dots, n\}),
\nonumber\end{align}
where
$v_k \in {\mathbb R}^{n}$
is the vector which has entries 1 in the
$k$th and $k+1$st positions (in positions $n$ and $1$ if $k=n$) and 0 otherwise.
Obviously,
$\S_n^T v_k = 0$
$(k \in \{1, \dots, n\})$
and this holds for all linear combinations of the
$v_k$,
too. In fact, the vectors
$\{v_1, v_2, \dots, v_{n-1}\}$
span the subspace
$\{\S_n\}^\bot$:
given
$u \in \{\S_n\}^\bot$,
we can take
$\alpha_1 = u_1$,
$\alpha_2 = u_2 - u_1$,
$\alpha_3 = u_3 - u_2 + u_1$,
etc., ending with
$\alpha_{n} = u_{n} - u_{n-1} + u_{n-2} - \cdots + u_2 - u_1 = 0;$
then
$u = \sum_{j=1}^{n-1} \alpha_j v_j.$
Therefore, by bilinearity a square matrix with property (M) also satisfies (f).
The converse is straightforward.

To see that (g) is equivalent to the condition (N), consider the vectors
$v_k$ defined in the proof of
%$v_k \in \mathbb{R}^n$, $k \in \{1, \dots, n-1\}$, defined in the proof of
part (f), which span the space $\{\S_n\}^\bot$.
As $M \S_n$ and $\S_n^T M$ are the vectors of alternating row and column sums of
$M$, respectively, (N) implies that
\begin{equation*}
 \S_n^T M v_k = 0 = v_k^T M \S_n
 \qquad (k \in \{1, \dots, n-1\}),
\end{equation*}
and hence (g) by linearity; the converse is obvious.
\phantom{.}\hfill$\qed$\par\noindent

\medskip\noindent
We now observe that the conditions in Theorem \ref{tsymeq} (a) and (e), as well as
those in (f) and (g), are essentially mutually exclusive.
\begin{Lemma}
\label{lMNSVdisj}
Let $n \in \mathbb{N}$ and $y \in \mathbb{R}^n\setminus\{0_n\}$.
If $M \in \mathbb{R}^{n \times n}$ satisfies

(i) $y^T M u = 0 = u^T M y$
$(u \in \{y\}^\bot)$,

(ii) $u^T M v = 0 $ $(u, v \in \{y\}^\bot)$, and

(iii) $y^T M y = 0$,

\noindent
then $M = \omx_n$.
\end{Lemma}
\noindent
{\it Proof.} The matrix
$P = (y^T y)^{-1}\,y y^T$ is symmetric and idempotent, $P^2 = P$; it follows
that $\imx_n - P$ also has these properties.
If $u \in \mathbb{R}^n$, then $P u$ is a multiple of $y$ and
$y^T P u = y^T u$, so $(\imx_n - P) u \in \{y\}^\bot$.
Hence, for $u, v \in \mathbb{R}^n$,
\begin{equation*}
 u^T M v = (P u)^T M P v + ((\imx_n - P) u)^T M P v
  + (P u)^T M (\imx_n - P) v  + ((\imx_n - P) u)^T M (\imx_n - P) v
 = 0,
\end{equation*}
where the first term vanishes by (iii), the second and third by (i) and the
fourth by (ii).
\phantom{.}\hfill$\qed$\par\noindent
\begin{Theorem}
\label{tdisjoint}
Let $n \in \mathbb{N}$. Then
$A_n \cap B_n = \{\omx_n\}$, $S_n \cap V_n = \{\omx_n\}$ and
$M_n \cap N_n = \{\omx_n\}$.
If $n$ is even, then also
$P_n \cap Q_n = \{\omx_n\}$.
\end{Theorem}
\noindent
{\it Proof.\/}
The second and third identity follow from Lemma \ref{lMNSVdisj}, taking
$y = 1_n$ and $y = \S_n$, respectively.
The first and fourth identity are immediate from combining (A), (B) and (P), (Q),
respectively, with weight $w = 0$.
\phantom{.}\hfill$\qed$\par\noindent

\begin{Corollary}
\label{cdirsum}
Let $n \in \mathbb{N}$. Then
$\mathbb{R}^{n\times n} = B_n \oplus A_n = S_n \oplus V_n
= N_n \oplus M_n$.

If $n$ is even, then also
$\mathbb{R}^{n\times n} = Q_n \oplus P_n$.
\end{Corollary}

\noindent
{\it Proof.\/}
In view of Theorem \ref{tdisjoint}, we only need to show that any $n\times n$
matrix can be written as a sum of matrices from each summand in all cases.

Let $M \in \mathbb{R}^{n\times n}$.
Then
\[
 M = \frac{1}{2} (M + \jmx_n M\jmx_n) + \frac{1}{2} (M - \jmx_n M\jmx_n),
\]
and using Theorem \ref{tsymeq} (b), (c) and the fact that $\jmx_n^2 = \imx_n$, we
see that the first term is in $B_n$, the second in $A_n$.

Further, defining the projector $P$ as in the proof of Lemma \ref{lMNSVdisj},
we find
\[
 M = (P M P + (\imx_n - P) M (\imx_n - P)) + (P M (\imx_n - P) + (\imx_n - P) M P);
\]
then for $y = 1_n$, the first bracket is in $S_n$, the second in $V_n$ by
Theorem \ref{tsymeq} (a), (e);
for $y = \S_n$, the first bracket is in $N_n$, the second in $M_n$ by Theorem
\ref{tsymeq} (g), (f).

Finally, if $n = 2 \nu$ is even, then we can split $M$ into $\nu\times\nu$ blocks,
\[
 M = \pmatrix{A & B \cr C & D} = \frac{1}{2} \pmatrix{A+D & B+C \cr B+C & A+D}
   + \frac{1}{2}\pmatrix{A-D & B-C \cr -(B-C) & -(A-D)},
\]
with the first matrix on the right-hand side, of form (\ref{eQstruct}), in $Q_n$,
the second matrix, of form (\ref{ePstruct}), in $P_n$.
\phantom{.}\hfill$\qed$\par\noindent

\section{Representation Formulae: the Type B+A and Q+P Algebras}
\label{sRFBQ}
We now proceed to finding representation formulae for the various symmetry types.
These will give a way of constructing matrices of a particular symmetry type by an
expression with no or much simpler constraints; in the cases of even-dimensional
type S, V, N or M matrices there will be a recursive element in that the construction
formula requires some lower-dimensional matrix of the same type.
Furthermore, these representation formulae will make the relationship between
symmetry types and their algebraic properties more transparent.

We start with types A and B.
As a template for this approach, consider the characterisation and construction
of combined type A, S and type B, S matrices in the paper \cite{hill}.

A crucial role is played by the matrix $\xmx_n$, which is used to transform square
matrices to their {\em block representation} by conjugation (see \cite{hill} and, more generally, \cite{njh}); it takes the form
\begin{equation*}
 \xmx_{n} = \rc{\sqrt 2} \pmatrix{ \imx_\nu & \jmx_\nu \cr \jmx_\nu & -\imx_\nu} \in {\mathbb R}^{n \times n}
\nonumber\end{equation*}
if $n = 2 \nu$ is even, and
\begin{equation*}
 \xmx_{n} = \pmatrix{\rc{\sqrt 2}\imx_\nu & 0_\nu & \rc{\sqrt 2}\jmx_\nu \cr 0_\nu^T & 1 & 0_\nu^T \cr
         \rc{\sqrt 2}\jmx_\nu & 0_\nu & -\rc{\sqrt 2}\imx_\nu}
 \in {\mathbb R}^{n \times n}
\nonumber\end{equation*}
if $n = 2 \nu+1$ is odd.
The matrix
$\xmx_{n}$
is an orthogonal symmetric involution, i.e.
$\xmx_{n}^T = \xmx_{n}$
and
$\xmx_{n}^2 = \imx_{n}.$
It follows that
\begin{equation}
 (\xmx_n M \xmx_n) (\xmx_n M' \xmx_n) = \xmx_n (M M') \xmx_n
 \qquad (M, M' \in \mathbb{R}^{n \times n}),
\end{equation}
so conjugation with $\xmx_n$ (which is also linear) is a matrix algebra
homomorphism.

\noindent

Specifically for the weight matrix $\emx_n$, the block representation is
\begin{equation*}
 \emx_{n} = \xmx_{n} \pmatrix{2 \emx_\nu & \omx_\nu \cr \omx_\nu & \omx_\nu} \xmx_{n}
\nonumber\end{equation*}
if $n = 2\nu$ is even, and
\begin{equation}
\label{eemxblodd}
 \emx_{n} = \xmx_{n} \pmatrix{2 \emx_\nu & \sqrt 2\,1_\nu & \omx_\nu \cr \sqrt 2\,1_\nu^T & 1 & 0_\nu^T
     \cr \omx_\nu & 0_\nu & \omx_\nu} \xmx_{n},
\end{equation}
if $n = 2\nu + 1$ is odd.

\begin{Lemma}\label{ttopsy}
A matrix
$M \in \mathbb{R}^{n \times n}$ is an element of $A_n$
if and only if
\begin{align}
 M &= \xmx_n \pmatrix{\omx & \Psi \cr \Phi & \omx_\nu} \xmx_n,
\label{ePQblock}\end{align}
where
$\Phi, \Psi \in \mathbb{R}^{\nu\times\nu}$
if $n = 2\nu$ is even,
$\Phi \in \mathbb{R}^{\nu\times(\nu+1)}$,
$\Psi \in \mathbb{R}^{(\nu+1)\times\nu}$
if $n = 2\nu+1$ is odd, and the top left null matrix has matching size.
\end{Lemma}
\par\medskip\noindent
{\it Proof. \/}
In the case of even $n$, it follows from Theorem \ref{tsymeq} (b) and
\begin{equation}
\label{eJeven}
 \jmx_n = \pmatrix{\omx_\nu & \jmx_\nu \cr \jmx_\nu & \omx_\nu}
\end{equation}
that we
can write the weight zero generally associated matrix in the form
\begin{align}
 M &= \pmatrix{A & -\jmx_\nu B \jmx_\nu \cr B & - \jmx_\nu A \jmx_\nu}
\nonumber\end{align}
with some $A, B \in \mathbb{R}^{\nu\times\nu}$;
then its block representation is
\begin{align}
 \xmx_n \pmatrix{A & -\jmx_\nu B \jmx_\nu \cr B & - \jmx_\nu A \jmx_\nu} \xmx_n &= \pmatrix{\omx_\nu & A \jmx_\nu + \jmx_\nu B \jmx_\nu \cr \jmx_\nu A - B & \omx_\nu}.
\nonumber\end{align}
Conversely,
\begin{align}
 \xmx_n \pmatrix{\omx_\nu & \Psi \cr \Phi & \omx_\nu} \xmx_n &= \rc 2 \pmatrix{\Psi \jmx_\nu + \jmx_\nu \Phi & -(\Psi \jmx_\nu - \jmx_\nu \Phi) \jmx_\nu \cr
        \jmx_\nu (\Psi \jmx_\nu - \jmx_\nu \Phi) & -\jmx_\nu (\Psi \jmx_\nu + \jmx_\nu \Phi) \jmx_\nu},
\nonumber\end{align}
which evidently satisfies Theorem \ref{tsymeq} (b).

In the case of odd $n$, we have
\begin{equation}
\label{eJodd}
 \jmx_n = \pmatrix{\omx_\nu & 0_\nu & \jmx_\nu \cr 0_n^T & 1 & 0_n^T \cr
     \jmx_\nu & 0_\nu & \omx_\nu}
\end{equation}
and the matrix is of the form
\begin{equation*}
 M = \pmatrix{A & v & - \jmx_\nu B \jmx_\nu \cr w^T & 0 & -w^T \jmx_\nu \cr
              B & -\jmx_\nu v & -\jmx_\nu A \jmx_\nu},
\end{equation*}
where $A, B \in \mathbb{R}^{\nu\times\nu}$ and $v, w \in \mathbb{R}^\nu$.
Then its block representation is
\begin{equation*}
 \xmx_n \pmatrix{A & v & - \jmx_\nu B \jmx_\nu \cr w^T & 0 & -w^T \jmx_\nu \cr
              B & -\jmx_\nu v & -\jmx_\nu A \jmx_\nu} \xmx_n
  = \pmatrix{ \omx_\nu & 0_\nu & A\jmx_\nu + \jmx_\nu B \jmx_n \cr
              0_\nu^T & 0 & \sqrt 2 w^T \jmx_\nu \cr
              \jmx_\nu A - B & \sqrt 2 \jmx_\nu v_\nu & \omx_\nu}.
\end{equation*}
For the converse, the relationship between $A$, $B$ and the first $\nu$ columns of
$\Phi$ and rows of $\Psi$ is as in the even-dimensional case.
\phantom{.}\hfill$\qed$\par\noindent
\begin{Lemma}\label{tbalsy}
A matrix
$M \in \mathbb{R}^{n \times n}$ is an element of $B_n$
if and only if
\begin{align}
 M &= \xmx_n \pmatrix{\Upsilon & \omx \cr \omx & \Omega} \xmx_n
\label{eYOblock}\end{align}
with matrices
$\Upsilon, \Omega \in \mathbb{R}^{\nu\times\nu}$
if $n = 2\nu$ is even,
$\Upsilon \in \mathbb{R}^{(\nu+1)\times(\nu+1)}$,
$\Omega\in \mathbb{R}^{\nu\times\nu}$
if $n = 2\nu+1$ is odd, and null matrices of matching size.
\end{Lemma}
\par\medskip\noindent
{\it Proof. \/}
In the case of even $n$, it follows from Theorem \ref{tsymeq} (c) and (\ref{eJeven}) that we
can write the matrix in the form
\begin{align}
 M &= \pmatrix{A & \jmx_\nu B \jmx_\nu \cr B &  \jmx_\nu A \jmx_\nu}
\nonumber\end{align}
with $A, B \in \mathbb{R}^{\nu\times\nu}$;
then its block representation is
\begin{align}
 \xmx_n \pmatrix{A & \jmx_\nu B \jmx_\nu \cr B &  \jmx_\nu A \jmx_\nu} \xmx_n
 &= \pmatrix{A + \jmx_\nu B & \omx_\nu \cr \omx_\nu & \jmx_\nu A \jmx_\nu - B \jmx_\nu}.
\nonumber\end{align}
Conversely,
\begin{align}
 \xmx_n \pmatrix{\Upsilon & \omx_\nu  \cr  \omx_\nu & \Omega} \xmx_n
 &= \rc 2 \pmatrix{
\Upsilon + \jmx_\nu \Omega \jmx_\nu & \jmx_\nu(\jmx_\nu \Upsilon - \Omega \jmx_\nu)\jmx_\nu \cr \jmx_\nu \Upsilon - \Omega \jmx_\nu & \jmx_\nu \Upsilon \jmx_\nu + \Omega}
\nonumber\end{align}
clearly gives a balanced matrix.

In the case of odd $n$, the matrix is of the form
\begin{equation*}
 M = \pmatrix{A & v &  \jmx_\nu B \jmx_\nu \cr w^T & x & w^T \jmx_\nu \cr
              B & \jmx_\nu v & \jmx_\nu A \jmx_\nu},
\end{equation*}
with $A, B \in \mathbb{R}^{\nu\times\nu}$ and $v, w \in \mathbb{R}^\nu$, $x \in \mathbb{R}$.
Then its block representation is
\begin{equation*}
 \xmx_n \pmatrix{A & v &  \jmx_\nu B \jmx_\nu \cr w^T & x & w^T \jmx_\nu \cr
              B & \jmx_\nu v & \jmx_\nu A \jmx_\nu} \xmx_n
  = \pmatrix{ A + \jmx_\nu B & \sqrt 2 v & \omx_\nu\cr
              \sqrt 2 w^T & x & 0_n^T \cr
              \omx_\nu & 0_n & \jmx_\nu A \jmx_\nu - B \jmx_\nu}.
\end{equation*}
For the converse, the relationship between $A$, $B$ on the one hand and the top $\nu\times\nu$ submatrix of
$\Upsilon$ and the matrix $\Omega$ is as in the even-dimensional case;
$v, w$ and $x$ can be read off directly.
\phantom{.}\hfill$\qed$

\bigskip\noindent
The block representations of Lemma \ref{ttopsy} and \ref{tbalsy} make the splitting of a
general matrix into types A and B very transparent. Moreover, it becomes obvious
that this splitting gives
$\mathbb{R}^{n\times n} = B_n \oplus A_n$
the structure of a $\mathbb{Z}_2$-graded algebra, with `even' subalgebra $B_n$ and `odd' complement $A_n$.

\begin{Theorem}
\label{tBAgrad}
Let $n \in \mathbb{N}$.
Then
\begin{equation*}
 B_n B_n \subset B_n, \qquad
 A_n A_n \subset B_n, \qquad
 A_n B_n \subset A_n, \qquad
 B_n A_n \subset A_n. \qquad
\end{equation*}
\end{Theorem}
\noindent
The same structure can be seen in the $Q_n$ and $P_n$ symmetry types; indeed, a
straightforward calculation using directly the structures (\ref{eQstruct}) and
(\ref{ePstruct}) shows that
$\mathbb{R}^{n\times n} = Q_n \oplus P_n$
also is a $\mathbb{Z}_2$-graded algebra, with `even' subalgebra $Q_n$.
\begin{Theorem}
\label{tQPgrad}
Let $n \in \mathbb{N}$.
Then
\begin{equation*}
 Q_n Q_n \subset Q_n, \qquad
 P_n P_n \subset Q_n, \qquad
 P_n Q_n \subset P_n, \qquad
 Q_n P_n \subset P_n. \qquad
\end{equation*}
\end{Theorem}

\section{Representation Formulae: the Type S+V and N+M Algebras}
Although the block representation by conjugation with the matrix $\xmx_n$ was originally
devised to capture the structure of matrices with (A) or (B) symmetry,
it also proves useful in the study of other symmetry types.
\begin{Theorem}
\label{tSblock}
If $n = 2\nu$ is even,
$M \in \mathbb{R}^{n \times n}$ is an element of $S_n$ if and only if
\begin{align}
\label{eSbleven}
 M &= \xmx_n \pmatrix{Y & V^T \cr W & Z} \xmx_n
\end{align}
with $Y \in S_\nu$,
$V, W \in \mathbb{R}^{\nu \times \nu}$ with row sums $0$, and
$Z \in \mathbb{R}^{\nu \times \nu}$.

If $n = 2\nu+1$ is odd, $M \in \mathbb{R}^{n\times n}$ is an element of $S_n$
if and only if
\begin{equation}
 M = \xmx_{n} \pmatrix{Y + 2 w \emx_\nu & \sqrt 2 (w 1_\nu - Y 1_\nu) & V^T \cr
   \sqrt 2 (w 1_\nu - Y^T 1_\nu)^T & w + 2\,1_\nu^T Y 1_\nu & - \sqrt 2 (V 1_\nu)^T \cr
   W & -\sqrt 2 W 1_\nu & Z} \xmx_{n}
\label{eSblodd}\end{equation}
with arbitrary
$V, W, Y, Z \in {\mathbb R}^{\nu \times \nu}$; $w \in \mathbb{R}$ is the weight.
\end{Theorem}

\noindent
{\it Proof.}
First, consider the case of even $n$.
Then
\begin{equation}
\label{eonebleven}
\xmx_n 1_n = \pmatrix{\sqrt{2}\,1_\nu \cr 0_\nu}
\end{equation}
and hence
$u\in \{1_n\}^\bot$ if and only if $\xmx_n u = \pmatrix{\xi\cr\eta}$ with
$\xi\in\{1_\nu\}^\bot$ and arbitrary $\eta\in\mathbb{R}^\nu$.
Writing the block representation of $M$ in the form (\ref{eSbleven}),
we find that the conditions of Theorem \ref{tsymeq} (a) take the form
\begin{align*}
 0 &= \pmatrix{\xi\cr\eta}^T \pmatrix{Y & V^T \cr W & Z} \pmatrix{1_\nu \cr 0_\nu}
 = \xi^T Y 1_\nu + \eta^T W 1_\nu, \\\\
 0 &= \pmatrix{1_\nu \cr 0_\nu}^T \pmatrix{Y & V^T \cr W & Z} \pmatrix{\xi\cr\eta}
 = 1^T Y \xi + (V 1_\nu)^T \eta
\end{align*}
for all $\xi\in\{1_\nu\}$ and for any $\eta\in\mathbb{R}^\nu$.
This is equivalent to $V 1_\nu = 0_\nu$, $W 1_\nu = 0_\nu$ and (again by Theorem \ref{tsymeq}
(a)) $Y\in S_\nu$.

The case of odd $n$ is a bit more tricky. By Lemma \ref{lMNSVeigen}, $M \in S_n$ is equivalent to
$1_n$ being an eigenvector, for eigenvalue $w$, of both $M$ and $M^T$.
Hence, considering the matrix $M_0 := M - \frac{w}{n}\emx_n$, we find
$M_0 1_n = 0_n$, $M_0^T 1_n = 0_n$.
Now observing that
\begin{equation}
\label{eoneblodd}
 \xmx_n 1_n = \pmatrix{\sqrt 2\,1_\nu \cr 1 \cr 0_\nu}
\end{equation}
and writing the block representation of $M_0$ in the form
\begin{equation}
\label{egenblodd}
 M_0 = \xmx_n \pmatrix{Y & v & V^T \cr y^T & \alpha & z^T \cr W & x & Z}\xmx_n
\end{equation}
with
$V, W, Y, Z \in \mathbb{R}^{\nu\times\nu}$,
$x, y, v, z \in \mathbb{R}^\nu$ and $\alpha\in\mathbb{R}$,
we see that these conditions on $M_0$ are equivalent to
\begin{equation*}
 \pmatrix{0_\nu \cr 0 \cr 0_\nu} = \pmatrix{\sqrt 2\,Y 1_\nu + v \cr
   \sqrt 2\,y^T 1_\nu + \alpha \cr \sqrt 2\,W 1_\nu + x}, \qquad
 \pmatrix{0_\nu \cr 0 \cr 0_\nu} = \pmatrix{\sqrt 2\,Y^T 1_\nu + y \cr
   \sqrt 2\,v^T 1_\nu + \alpha \cr \sqrt 2\,V 1_\nu + z}.
\end{equation*}
Thus $x = -\sqrt 2\,W 1_\nu$, $y = -\sqrt 2\,Y^T 1_\nu$,
$v = -\sqrt 2\,Y 1_\nu$, $z = -\sqrt 2\, V 1_\nu$ and
$\alpha = -\sqrt 2\,1_\nu^T v = 2\,1_\nu^T Y 1_\nu = -\sqrt 2\,y^T 1_\nu$.
This gives (\ref{eSblodd}) in view of the block representation of $\emx_n$,
eq. (\ref{eemxblodd}).
\phantom{.}\hfill$\qed$
\begin{Theorem}
\label{tVblock}
If $n = 2 \nu$ is even, then $M \in \mathbb{R}^{n \times n}$ is an element of $V_n$
if and only if
\begin{equation}
\label{eVbleven}
M = \xmx_n \pmatrix{Y & 1_\nu a^T \cr b 1_\nu^T & \omx_\nu} \xmx_n
\end{equation}
with $Y \in V_\nu$ and $a, b \in \mathbb{R}^\nu$.

If $n = 2 \nu + 1$ is odd, then $M \in \mathbb{R}^{n \times n}$ is an element of
$V_n$ if and only if
\begin{equation}
\label{eVblodd}
M = \xmx_n \pmatrix{\sqrt 2 (v 1_\nu^T + 1_\nu y^T) - \frac{2\sqrt 2}{2\nu - 1}\,(1_\nu^T (v+y))\emx_\nu & v & \sqrt 2\, 1_\nu z^T \cr
y^T & \frac {\sqrt 2} {2\nu - 1}\,1_\nu^T(v+y) & z^T\cr
\sqrt 2\, x 1_\nu^T & x & \omx_\nu} \xmx_n
\end{equation}
with arbitrary $v, x, y, z \in \mathbb{R}^\nu$.
\end{Theorem}

\noindent
{\it Remark.\/}
As the right-most $\nu$ columns (and the bottom $\nu$ rows) of the block
representation form a matrix of rank at most 1, it follows that the rank of a
type V matrix cannot exceed $n - \nu + 1$. In fact, the rank cannot exceed 7, as it has a contribution of at most 3 from the top left block, 1 from the bottom left block, 1 from the middle row and also from the middle column, and 1 from the right $\nu$ columns. The combined contribution thus gives a maximum possible rank of 7.

We note that for both the type S and the type V matrices, that while the even-dimension formula looks simpler, it has a
recursive condition on $Y$, whereas the odd-dimension formula has no restrictions.
We also note that if $n$ is odd, then
the dimension of $V_n$ is $4\nu = 2 n - 2$. If $n$ is even, then the dimension of $V_n$
is the dimension of $V_{n/2}$ plus $n$; as $V_2$ has dimension 2, this also
works out as $2 n - 2$, which correctly implies that the dimension of $S_n$ is given by $n^2 - 2 n + 2$ (see also \cite{brunnock}).

\medskip
\noindent
{\it Proof.\/}
For even $n$, the first condition in Theorem \ref{tsymeq} (e) translates, in analogy
to the beginning of the proof of Theorem \ref{tSblock}, into
\begin{equation*}
 0 = \pmatrix{\xi\cr\eta}^T \pmatrix{Y & V^T \cr W & Z} \pmatrix{\xi'\cr\eta'}
 = \xi^T Y \xi' + \xi^T V^T \eta' + \eta^T W \xi' + \eta^T Z \eta',
\end{equation*}
for any $\xi,\xi' \in \{1_\nu\}^\bot$ and $\eta,\eta'\in\mathbb{R}^\nu$, in the block
representation.
When we take $\xi = \xi' = 0_\nu$, this implies $Z = \omx_\nu$. Taking one or both
of $\eta$, $\eta'$ to be $0_\nu$, we find
$V \xi = 0$, $W \xi = 0$ $(\xi \in \{1_\nu\}^\bot)$
and $\xi^T Y \xi' = 0$ $(\xi, \xi' \in \{1_\nu\}^\bot)$.
Similarly, the second condition in Theorem \ref{tsymeq} (e) gives
$1_\nu^T Y 1_\nu = 0$.
Hence $Y \in V_\nu$, and as the rows of $V, W$ must be elements of
$\{1_\nu\}^{\bot\bot} = \{1_\nu\}$, these matrices are of the stated form.

For odd $n$, the first condition of Theorem \ref{tsymeq} (e) takes the block form
\begin{align}
 0 &= \pmatrix{\xi \cr - \sqrt 2\, 1_\nu^T \xi \cr \eta}
     \pmatrix{Y & v & V^T \cr y^T & \alpha & z^T \cr W & x & Z}
    \pmatrix{\xi' \cr - \sqrt 2\, 1_\nu^T \xi' \cr \eta'} \nonumber\\
 &= \xi^T Y \xi' - \sqrt 2\, (\xi^T v)(1_\nu^T \xi') + \xi V^T \eta'
  - \sqrt 2\,(\xi^T 1_\nu)(y^T \xi') + 2 (\xi^T 1_\nu)(1\_\nu^T \xi')\alpha
\nonumber\\
 &\qquad - \sqrt 2\,(\xi^T 1_\nu)(z^T \eta') + \eta^T W \xi' - \sqrt 2\,(\eta^T x)(1_\nu^T\xi')
 + \eta^T Z \eta',
\label{eVa}
\end{align}
for all $\xi, \xi', \eta, \eta' \in \mathbb{R}^\nu$.

Taking $\xi = \xi' = 0_\nu$, we conclude that $Z = \omx_\nu$.
Then, taking $\xi = 0_\nu$, we see that
$\eta^T W \xi' = \sqrt 2\, \eta^T x 1_\nu^T \xi'$ for all $\eta, \xi' \in\mathbb{R}^\nu$, which implies $W = \sqrt 2\, x 1^T$.
Similarly, taking $\xi' = 0_\nu$ gives $V = \sqrt 2\, z 1^T$.
This leaves (\ref{eVa}) in the form
\begin{equation*}
 0 = \xi^T Y \xi' - \xi^T (\sqrt 2\, v 1_\nu^T + \sqrt 2\, 1_\nu y^T - 2 \alpha
      1_\nu 1_\nu^T) \xi'
\end{equation*}
for all $\xi, \xi' \in \mathbb{R}^\nu$, so
$Y = \sqrt 2 (v 1_\nu^T + 1_\nu y^T) - 2 \alpha \emx_\nu$.

Furthermore, the second condition in Theorem \ref{tsymeq} (e) takes the block
representation form
\begin{align*}
 0 &= \pmatrix{\sqrt 2\,1_\nu \cr 1 \cr 0_\nu}
     \pmatrix{Y & v & V^T \cr y^T & \alpha & z^T \cr W & x & Z}
     \pmatrix{\sqrt 2\,1_\nu \cr 1 \cr 0_\nu}
 = 2\,1_\nu^T Y 1_\nu + \sqrt 2\, 1_\nu^T v + \sqrt 2\, y^T 1_\nu + \alpha,
\end{align*}
which together with the previous identity for $Y$ gives
$\alpha = \frac{\sqrt 2}{2\nu - 1}\,(1_\nu^T v + y^T 1)$, and hence (\ref{eVblodd}).

It is a straightforward calculation to check that, conversely, (\ref{eVblodd})
satisfies Theorem \ref{tsymeq} (e).
\phantom{.}\hfill$\qed$

\bigskip\noindent
Similar considerations yield the block representation, and hence a construction
method, for type N and M matrices.
\begin{Theorem}
\label{tNblock}
If $n = 2\nu$ is even, then $M \in \mathbb{R}^{n \times n}$ is an element of
$N_n$ if and only if
\begin{equation}
\label{eNbleven}
 M = \xmx_n \pmatrix{Y & V^T \cr W & Z} \xmx_n
\end{equation}
with $Y\in\mathbb{R}^{\nu\times\nu}$, $V, W\in\mathbb{R}^{\nu\times\nu}$ such that
$V^T\S_\nu = W^T\S_\nu = 0_\nu$, and $Z \in N_\nu$.

If $n = 2\nu+1$ is odd, then $M\in \mathbb{R}^{n \times n}$ is an element of
$N_n$ if and only if
\begin{equation}
\label{eNblodd}
 M = \xmx_n \pmatrix{Y + 2\lambda \S_\nu\S_\nu^T & \pm\sqrt 2\,(\lambda\S_\nu - Y\S_\nu) & V^T \cr
 \pm \sqrt 2\,(\lambda\S_\nu - Y^T\S_\nu)^T & \lambda + 2\,\S_\nu^T Y \S_\nu & \mp\sqrt 2\,(V\S_\nu)^T \cr
 W & \mp\sqrt 2\,W\S_\nu & Z} \xmx_n
\end{equation}
with arbitrary $V, W, Y, Z \in\mathbb{R}^{\nu\times\nu}$ and $\lambda\in\mathbb{R}$;
here the upper signs apply if $\nu$ is even, the lower signs if $\nu$ is odd.

It follows that in the even-dimensional case, $V$ and $W$ have the rather strange property of {\it column\/}
(alternating) sum 0, which is different from type S and also from the
odd-dimensional case, where there is no restriction on column sums and
an effective weighted row sum 0.
\end{Theorem}

\medskip\noindent
{\it Proof.\/} The proof of Theorem \ref{tNblock} is largely analogous to the proof
of Theorem \ref{tSblock}, with the vector $\S_n$ taking the role of the vector $1_n$,
so we just detail the differences.

In the case of even $n$, we note that $J_\nu \S_\nu = \mp \S_\nu$ (where the
upper sign applies if $\nu$ is even, the lower sign if $\nu$ is odd), so
\begin{equation}
\label{eSSbleven}
 \xmx_n \S_n = \frac 1 {\sqrt 2} \pmatrix{\imx_\nu & \jmx_\nu \cr \jmx_\nu & -\imx_\nu} \pmatrix{\S_\nu \cr \pm \S_\nu}
 = \mp \sqrt 2\pmatrix{0_\nu \cr \S_\nu}.
\end{equation}
Thus $u \in \{\S_n\}^\bot$ if and only if $\xmx_n u = \pmatrix{\xi\cr\eta}$
with arbitrary $\xi\in \mathbb{R}^\nu$ and $\eta\in\{\S_\nu\}^\bot$.
The conditions of Theorem \ref{tsymeq} (g) take the form
\begin{equation*}
0 = \xi^T V^T\S_\nu + \eta^T Z \S_\nu, \qquad 0 = \S_\nu^T W \xi + \S_\nu^T Z\eta,
\end{equation*}
and give the conditions on $V, W, Z$ stated in the theorem.

In the case of odd $n$,
we apply Lemma \ref{lMNSVeigen} with $y = \S_n$ to find that $\S_n$ is an
eigenvector, for eigenvalue $\lambda\in\mathbb{R}$, of both $M$ and $M^T$,
and consider the matrix $M_0 := M - \frac \lambda n\,\S_n\S_n^T$. Then
$M_0 \S_n = 0_n$, $M_0 \S_n = 0_n$.
Now
\begin{equation}
\label{eSSblodd}
 \xmx_n \S_n = \xmx_n \pmatrix{\S_\nu \cr \pm 1 \cr \mp \S_\nu}
 = \pmatrix{\sqrt 2\,\S_\nu \cr \pm 1 \cr 0_n},
\end{equation}
and writing the block representation of $M_0$ as in (\ref{egenblodd}), we see
that the conditions on $M_0$ are equivalent to
\begin{equation*}
 \pmatrix{0_\nu \cr 0 \cr 0_\nu} = \pmatrix{\sqrt 2\,Y\S_\nu \pm v \cr
  \sqrt 2\, y^T\S_\nu \pm \alpha \cr \sqrt 2\, W\S_\nu \pm x},
 \qquad
 \pmatrix{0_\nu \cr 0 \cr 0_\nu} = \pmatrix{\sqrt 2\,Y^T\S_\nu \pm y \cr
  \sqrt 2\, v^T\S_\nu \pm \alpha \cr \sqrt 2\, V\S_\nu \pm z},
\end{equation*}
from which $v, x, y, z$ and $\alpha$ can be expressed in terms of $V, W$ and $Y$.
Equation (\ref{eNblodd}) follows by observing that
\begin{equation*}
 \xmx_n \S_n \S_n^T \xmx_n = \xmx_n \S_n (\xmx_n \S_n)^T
 = \pmatrix{2\,\S_\nu\S_\nu^T & \pm\sqrt 2\,\S_\nu & 0_\nu \cr
  \pm\sqrt 2\,\S_\nu^T & 1 & 0_\nu^T \cr \omx_\nu & 0_\nu & \omx_\nu}.
\end{equation*}
\phantom{.}\hfill$\qed$

\begin{Theorem}
\label{tMblock}
If $n = 2\nu$ is even, then $M \in \mathbb{R}^{n\times n}$ is an element of
$M_n$ if and only if
\begin{equation}
\label{eMbleven}
 M = \xmx_n \pmatrix{\omx_\nu & a \S_\nu^T \cr \S_\nu b^T & Z} \xmx_n
\end{equation}
with $Z\in M_\nu$ and $a, b\in\mathbb{R}^\nu$.

If $n = 2\nu+1$ is odd, then $M \in \mathbb{R}^{n\times n}$ is an element of $M_n$
if and only if
\begin{equation}
\label{eMblodd}
M = \xmx_n \pmatrix{\pm\sqrt 2\,(v \S_\nu^T + \S_\nu y^T)\mp\frac{2\sqrt 2}{2\nu-1}
 \,(\S_\nu^T(v+y))\S_\nu\S_\nu^T & v & \pm\sqrt 2\,\S_\nu z^T \cr
 y^T & \frac{\pm\sqrt 2}{2\nu-1}\,\S_\nu^T(v+y) & z^T \cr
 \pm\sqrt 2\,x \S_\nu^T & x & \omx_\nu} \xmx_n
\end{equation}
with arbitrary $v, x, y, z \in\mathbb{R}^\nu$; the upper sign applies if $\nu$ is
even, the lower sign if $\nu$ is odd.

Hence the rank 1 matrices are the `other way around' compared to the type V rank 1 matrices, where
again the odd-dimensional case has a very different structure to that of the even-dimensional case.
\end{Theorem}

\medskip\noindent
{\it Proof.\/}
In the case of even $n$, we use the formal block representation (\ref{eNbleven}),
and by (\ref{eSSbleven}), the conditions of Theorem \ref{tsymeq} (f) become
\begin{align*}
 0 &= \xi^T Y \xi' + \xi^T V^T \eta' + \eta^T W \xi' + \eta^T Z \eta', \\
 0 &= \S_\nu^T Z \S_\nu,
\end{align*}
for all $\xi, \xi'\in\mathbb{R}^\nu$ and $\eta, \eta' \in \{\S_\nu\}^\bot$, so the
stated properties of $V, W, Y, Z$ follow as in the proof of Theorem \ref{tVblock}.

For odd $n$, the reasoning is very similar to the proof of Theorem \ref{tVblock}.
However, in view of (\ref{eSSblodd}) we now have $u\in\{\S_n\}^\bot$ if and only if
\begin{equation*}
\xmx_n u = \pmatrix{\xi\cr\mp\sqrt 2\,\S_\nu^T\xi\cr\eta},
\quad (\xi,\eta\in\mathbb{R}^\nu),
\end{equation*}
so the analogue of (\ref{eVa}) takes the form
\begin{align*}
 0 &= \xi^T Y \xi' \pm \sqrt 2\, (\xi^T v)(\S_\nu^T \xi') + \xi V^T \eta'
  \pm \sqrt 2\,(\xi^T \S_\nu)(y^T \xi') + 2 (\xi^T \S_\nu)(\S\_\nu^T \xi')\alpha
\nonumber\\
 &\qquad \pm \sqrt 2\,(\xi^T \S_\nu)(z^T \eta') + \eta^T W \xi' \pm \sqrt 2\,(\eta^T x)(\S_\nu^T\xi')
 + \eta^T Z \eta'
\end{align*}
for all $\xi,\eta\in\mathbb{R}^\nu$; the remaining calculations are as before.
\phantom{.}\hfill$\qed$

\medskip\noindent
The splittings $\mathbb{R}^{n\times n} = S_n \oplus V_n$ and
$\mathbb{R}^{n\times n} = N_n \oplus M_n$ again have the structure of
a $\mathbb{Z}_2$-graded algebra, with subalgebra $S_n$ and $N_n$, respectively.
This follows from the next, more general, observation, which uses the symmetry
properties in their matrix algebra form (Theorem \ref{tsymeq}) directly,
rather than the block representations of Theorems \ref{tSblock} and \ref{tVblock}.
\begin{Lemma}
\label{lMNSValg}
Let $n \in \mathbb{N}$ and $y \in \mathbb{R}^N\setminus\{0_n\}$.
Let (i), (ii), (iii) denote the conditions listed in Lemma \ref{lMNSVdisj}.
\begin{description}
\item{(a)} If $M, M' \in \mathbb{R}^{n\times n}$ either both satisfy condition (i)
or both satisfy conditions (ii) and (iii),
then $M M'$ satisfies condition (i).
\item{(b)}
If $M\in \mathbb{R}^{n\times n}$ satisfies condition (i) and
$M'\in \mathbb{R}^{n\times n}$ satisfies conditions (ii) and (iii), then
$M M'$ and $M' M$ satisfy conditions (ii) and (iii).
\end{description}
\end{Lemma}
\medskip\noindent
{\it Proof.\/}
Let $P$ be the projector defined in the proof of Lemma \ref{lMNSVdisj}, and
$u, v \in \{y\}^\bot$.
Then, for (a) we observe
\begin{equation*}
 y^T M M' u = y^T M P M' u + y^T M (\imx_n - P) M' u = 0,
\end{equation*}
as in each of the two situations one half of each term vanishes; and similarly
$u M M' y^T = 0$.
For (b), we note that
\begin{equation*}
 u^T M M' v = u^T M P M' v + u^T M (\imx_n - P) M' v = 0,
\end{equation*}
since $u^T M P = 0_n$ and $(\imx_n - P) M' v = 0_n$; and
\begin{equation*}
 y^T M M' y = y^T M P M' y + y^T M (\imx_n - P) M' y = 0,
\end{equation*}
since $P M' y = 0_n$ and $y^T M (\imx_n - P) = 0_n$.
\phantom{.}\hfill$\qed$

\medskip\noindent
In view of Theorem \ref{tsymeq} (a) and (e), the choice $y = 1_n$ immediately gives
the following result.
\begin{Theorem}
\label{tSValg}
Let $n \in \mathbb{N}$. Then
\begin{equation*}
S_n S_n \subset S_n, \qquad
V_n V_n \subset S_n, \qquad
V_n S_n \subset V_n, \qquad
S_n V_n \subset V_n. \qquad
\end{equation*}
\end{Theorem}

\noindent
Similarly, by Theorem \ref{tsymeq} (f) and (g), the choice $y = \S_n$ gives the
following statement.
\begin{Theorem}
\label{tNMalg}
Let $n \in \mathbb{N}$. Then
\begin{equation*}
N_n N_n \subset N_n, \qquad
M_n M_n \subset N_n, \qquad
M_n N_n \subset M_n, \qquad
N_n M_n \subset M_n. \qquad
\end{equation*}
\end{Theorem}

\section{Representation Formulae: the R Algebra}

We now turn to symmetry type (R). This does not directly fit into the scheme of
pairings we observed in the other symmetry types; nevertheless, the block
representation turns out to be a valuable tool for constructing type R matrices
and for understanding their properties.

\begin{Theorem}
\label{tRblock}
If $n = 2\nu$ is even, then $M\in\mathbb{R}^{n\times n}$ is an element of $R_n$
if and only if
\begin{equation}
\label{eRbleven}
M = \xmx_n \pmatrix{\gamma\emx_\nu & 1_\nu z^T \cr x 1_\nu^T & Z} \xmx_n
\end{equation}
with $Z \in \mathbb{R}^{\nu\times\nu}$, $x, z \in\mathbb{R}^\nu$ and $\gamma\in\mathbb{R}$.

If $n = 2\nu+1$ is odd, then $M\in\mathbb{R}^{n\times n}$ is an element of $R_n$
if and only if
\begin{equation}
\label{eRblodd}
M = \xmx_n \pmatrix{\sqrt 2\,\gamma\emx_\nu & \gamma 1_\nu & \sqrt 2\,1_\nu z^T \cr
    \gamma 1_\nu^T & \frac\gamma{\sqrt 2} & z^T \cr
    \sqrt 2\,x 1_\nu^T & x & Z} \xmx_n
\end{equation}
with $Z \in \mathbb{R}^{\nu\times\nu}$, $x, z \in\mathbb{R}^\nu$ and $\gamma\in\mathbb{R}$.
\end{Theorem}

\noindent
{\it Proof.\/}
By Theorem \ref{tsymeq} (d), $M\in R_n$ if and only if
$(M + \jmx_n M) \mathbb{R}^n \subset \mathbb{R} 1_n$ and
$(M^T + \jmx_n M^T) \mathbb{R}^n \subset \mathbb{R} 1_n$.

In the case of even $n$, this means, using the formal block representation
(\ref{eSbleven}), equation (\ref{eonebleven}) and
\begin{equation}
\label{eJbleven}
 \xmx_n \jmx_n \xmx_n = \pmatrix{\imx_\nu & \omx_\nu \cr \omx_\nu & -\imx_\nu},
\end{equation}
that
\begin{align*}
2 \pmatrix{Y & V^T \cr \omx_\nu & \omx_\nu} \mathbb{R}^n &=
\left ( \pmatrix{Y & V^T \cr W & Z}
 + \pmatrix{\imx_\nu & \omx_\nu \cr \omx_\nu & -\imx_\nu}\pmatrix{Y & V^T \cr W & Z}
\right) \mathbb{R}^n \subset \mathbb{R} \pmatrix{1_\nu \cr 0_\nu}, \\
2 \pmatrix{Y^T & W^T \cr \omx_\nu & \omx_\nu} \mathbb{R}^n &=
\left ( \pmatrix{Y^T & W^T \cr V & Z^T}
 + \pmatrix{\imx_\nu & \omx_\nu \cr \omx_\nu & -\imx_\nu}\pmatrix{Y^T & W^T \cr V & Z^T}
\right) \mathbb{R}^n \subset \mathbb{R} \pmatrix{1_\nu \cr 0_\nu},
\end{align*}
equivalent to all columns of $Y, Y^T, V^T$ and $W^T$ being multiples of $1_\nu$.

Similarly, in the case of odd $n$, we use (\ref{eoneblodd})
and the formal block representation (\ref{egenblodd}) along with
\begin{equation*}
\xmx_n \jmx_n \xmx_j = \pmatrix{\imx_\nu & 0_\nu & \omx_\nu \cr 0_\nu^T & 1 & 0_\nu^T
  \cr \omx_\nu & 0_\nu & -\imx_\nu}
\end{equation*}
to rewrite the above conditions as
\begin{align*}
2\pmatrix{Y & v & V^T \cr y^T & \alpha & z^T \cr \omx_\nu & 0_\nu & \omx_\nu}
 \mathbb{R} \subset \mathbb{R}\pmatrix{\sqrt 2\,1_\nu \cr 1 \cr 0_\nu}, &\qquad
2\pmatrix{Y^T & y & W^T \cr v^T &\alpha & x^T \cr \omx_\nu & 0_\nu & \omx_\nu}
 \mathbb{R} \subset \mathbb{R}\pmatrix{\sqrt 2\,1_\nu \cr 1 \cr 0_\nu}.
\end{align*}
Hence we conclude that $Y = \sqrt 2\,\gamma 1_\nu 1_\nu^T$ for some
$\gamma\in\mathbb{R}$, $v = y = \gamma 1_\nu$ and $\alpha = \frac\gamma{\sqrt 2}$.
Moreover, $V^T = \sqrt 2\,1_\nu z^T$ and $W^T = \sqrt 2\,1_\nu x^T$.
\phantom{.}\hfill$\qed$

\medskip
\noindent
From the block representations of Theorem \ref{tRblock}, it is apparent that
the type R matrices form a subalgebra of the matrix algebra $\mathbb{R}^{n\times n}$.

\begin{Theorem}
\label{tRalg}
Let $n \in \mathbb{N}$. Then $R_n R_n \subset R_n$.
\end{Theorem}

\noindent
{\it Proof.\/} It is sufficient to show that the product of block representations
of type R matrices is the block representation of a type R matrix.

Let $\gamma, \gamma'\in\mathbb{R}$, $x, z, x', z'\in\mathbb{R}^\nu$
and $Z, Z'\in\mathbb{R}^{\nu\times\nu}$.
Then, for even $n = 2\nu$,
\begin{equation*}
\pmatrix{\gamma\emx_\nu & 1_\nu z^T \cr b 1_\nu^T & Z}
\pmatrix{\gamma'\emx_\nu & 1_\nu z'^T \cr x' 1_\nu^T & Z'}
= \pmatrix{(\gamma\gamma'\nu + z^T x')\emx_\nu & 1_\nu (\gamma\nu z'^T + z^T Z') \cr
 (\gamma'\nu x + Z x') 1_\nu^T & \nu x z'^T + Z Z'},
\end{equation*}
which is of the form (\ref{eRbleven}).
For odd $n = 2\nu+1$,
\begin{align*}
\pmatrix{\sqrt 2\,\gamma\emx_\nu & \gamma 1_\nu & \sqrt 2\,1_\nu z^T \cr
    \gamma 1_\nu^T & \frac\gamma{\sqrt 2} & z^T \cr \sqrt 2\,x 1_\nu^T & x & Z}
&\pmatrix{\sqrt 2\,\gamma'\emx_\nu & \gamma' 1_\nu & \sqrt 2\,1_\nu z'^T \cr
    \gamma' 1_\nu^T & \frac{\gamma'}{\sqrt 2} & z'^T \cr \sqrt 2\,x' 1_\nu^T & x' & Z'}
\\
&=\pmatrix{(n\gamma\gamma'+2 z^T x')\emx_\nu & \frac{n\gamma\gamma'+2 z^T x'}{\sqrt 2}\, 1_\nu & 1_\nu (n \gamma z'^T + \sqrt 2\, z Z') \cr
\frac{n\gamma\gamma'+2 z^T x'}{\sqrt 2}\,1_\nu^T & \frac{n\gamma\gamma'+2 z^T x'}2
 & \frac{n \gamma z'^T + \sqrt 2\, z Z'}{\sqrt 2} \cr
(n\gamma' x + \sqrt 2\,Z x') 1_\nu^T & \frac{n\gamma' x + \sqrt 2\,Z x'}{\sqrt 2} &
 n x z'^T + Z Z'},
\end{align*}
which is of the form (\ref{eRblodd}).
\phantom{.}\hfill$\qed$

\medskip\noindent
\begin{rmk}[concerning the complementary space to $R_n$] Given the $\mathbb{Z}_2$ graded algebras formed from the symmetries of the direct sum pairs
$B_n \oplus A_n$, $N_n \oplus M_n$, $S_n \oplus V_n$, and $Q_n \oplus P_n$, it seems sensible to ask the questions `what is
the complementary space to $R_n$?' and `do they form a $\mathbb{Z}_2$-graded
algebra?'
The complement to the (even-dimension) block representation would be
\begin{equation*}
\pmatrix{Y & V^T \cr W & Z}
\end{equation*}
with $V 1_\nu = W 1_\nu = 0_\nu$, $1_\nu^T Y 1_\nu = 0$; or equivalently
$1_\nu^T Y + u^T W 1_\nu + 1_\nu^T V^T v (+ u^T Z v) = 0$ $(u, v\in\mathbb{R}^\nu)$.
As $\xmx_n \pmatrix{1_\nu \cr u} = 1_n + \pmatrix{\jmx_\nu u \cr -u}$,
this is equivalent to the matrix $M$ having the property that
$(1_n + u)^T M (1_n + v) = 0$ for all $u, v\in\mathbb{R}^n$ such that
$\jmx_n u = -u, \jmx_n v = -v$.
This defines another space which directly sums with $R_n$ to the whole matrix
space; however, on first appearance it does not seem to form a $\mathbb{Z}_2$-graded
algebra with $R_n$, although their relationship may well merit further examination.
\end{rmk}

\section{Composite Symmetry: Most Perfect Squares}
\label{sCSMPS}

\noindent
After studying the basic symmetry types defined in Section \ref{sMSTS}, we now
proceed to the more complicated symmetries of most perfect square matrices and,
in the next section,
reversible square matrices. As the spaces of these matrices arise as intersections
of some basic symmetry spaces, their algebraic properties as well as construction
formulae can be readily deduced from the results in the preceding sections.

We use the convention of calling a direct sum $\Xi\oplus{\mathrm H}$, where $\Xi, {\mathrm H}$ are vector
subspaces of $\mathbb{R}^{n\times n}$, a $\mathbb{Z}_2$-graded algebra if
the first direct summand
$\Xi$ is the `even' subalgebra and
the second direct summand
${\mathrm H}$ is the `odd' complement, i.e. if
\begin{equation*}
 \Xi \Xi \subset \Xi, \qquad
 {\mathrm H} \Xi \subset {\mathrm H}, \qquad
 \Xi {\mathrm H} \subset {\mathrm H}, \qquad
 {\mathrm H} {\mathrm H} \subset \Xi.
\end{equation*}
The following statements follow immediately from this definition.
%[Remark on triple product in H?]

\begin{Lemma}
\label{lalginters}
Let $n\in\mathbb{N}$.
\begin{description}
\item{(a)} If $\Xi\oplus{\mathrm H}\in\mathbb{R}^{n\times n}$ is a $\mathbb{Z}_2$-graded algebra and
$\Gamma \subset \mathbb{R}^{n\times n}$ is a matrix algebra, then
$(\Xi \cap \Gamma)\oplus({\mathrm H}\cap \Gamma)$ is a $\mathbb{Z}_2$-graded
algebra.
\item{(b)} If $\Xi\oplus{\mathrm H}, \Xi'\oplus{\mathrm H}' \subset\mathbb{R}^{n\times n}$
are $\mathbb{Z}_2$-graded algebras, then
$(\Xi \cap \Xi')\oplus({\mathrm H}\cap{\mathrm H}')$ is a $\mathbb{Z}_2$-graded
algebra.
\end{description}
\end{Lemma}

\noindent
We now consider the set of most perfect square matrices,
$MPS_n = M_n \cap P_n \cap S_n$.
We recall that this is indeed the space of all weightless most perfect square
matrices, although the second part of property (M), corresponding to the last
condition in Theorem \ref{tsymeq} (f), was not stipulated in the original
definition of most perfect squares; indeed it is implied by condition (P)
as follows.
Using the fact that
$\S_n = \pmatrix{\S_\nu \cr \pm\S_\nu}$,
where the upper sign always refers to the case of even $\nu$, the lower to the case
of odd $\nu$,
equation (\ref{ePstruct}) gives
\begin{equation*}
 \S_n^T M \S_n = \pmatrix{\S_\nu \cr \pm\S_\nu}^T \pmatrix{A & B \cr -B & -A}
  \pmatrix{\S_\nu \cr \pm\S_\nu}
 = 0.
%\label{eSMS}
\end{equation*}
Since $\emx_n$ is a most perfect square matrix, the general most perfect square
matrices form the space $MSP_n \oplus \mathbb{R}\emx_n$.

The elements of $MSP_n$ have the following block representation.
\begin{Theorem}
\label{tMPSblock}
Let
$M\in\mathbb{R}^{n\times n}$,
$n = 2 \nu$ even.
Then
$M\in MPS_n$
if and only if there are
vectors
$a, b \in \{1_\nu\}^\bot$
with
$ \jmx_\nu a = \mp a$, $\jmx_\nu b = \mp b$,
where the upper sign applies if $\nu$ is even, the lower sign if $\nu$ is odd,
and a matrix
$Z \in A_\nu \cap M_\nu$,
such that
\begin{align}
 M &= \xmx_n \pmatrix{\omx_\nu & a \S_\nu^T \cr \S_\nu b^T & Z \cr} \xmx_n.
\nonumber\end{align}
\end{Theorem}
\par\medskip\noindent
{\it Proof. \/}
Let
$M\in MPS_n$.
Combining the block representations of Theorems \ref{tSblock} and \ref{tMblock},
we find that
\begin{equation*}
 M = \xmx_n \pmatrix{\omx_\nu & a \S_\nu^T \cr \S_\nu b^T & Z \cr} \xmx_n,
\end{equation*}
where
$a^T 1_\nu = 0 = b^T 1_\nu$
and $Z \in M_\nu$.
From (\ref{ePstruct}), we see that
\begin{align}
 \xmx_n M \xmx_n &= \rc 2 \pmatrix{A + B \jmx_\nu - \jmx_\nu B - \jmx_\nu A \jmx_\nu & A \jmx_\nu - B - \jmx_\nu B \jmx_\nu + \jmx_\nu A \cr
                       \jmx_\nu A + \jmx_\nu B \jmx_\nu + B + A \jmx_\nu & \jmx_\nu A \jmx_\nu - \jmx_\nu B + B \jmx_\nu - A \cr},
\nonumber\end{align}
and the calculation
\begin{align}
 \jmx_\nu (\jmx_\nu A \pm B \pm \jmx_\nu B \jmx_\nu + A \jmx_\nu) \jmx_\nu &= A \jmx_\nu \pm \jmx_\nu B \jmx_\nu \pm B + \jmx_\nu A
\nonumber\end{align}
shows that
$\S_\nu a^T, \S_\nu b^T \in B_\nu$.
Thus, by Theorem \ref{tsymeq} (c),
$\S_\nu a^T = \jmx_\nu \S_\nu a^T \jmx_\nu = \mp \S_\nu a^T \jmx_\nu$,
and hence $a = \mp \jmx_\nu a$; and similarly for $b$.
Also,
\begin{align}
 \jmx_\nu Z \jmx_\nu &= \rc 2\, \jmx_\nu (\jmx_\nu A \jmx_\nu - \jmx_\nu B + B \jmx_\nu - A) \jmx_\nu = \rc 2\, A - B \jmx_\nu + \jmx_\nu B - \jmx_\nu A \jmx_\nu = -Z,
\nonumber\end{align}
so
$Z \in A_\nu$
by Theorem \ref{tsymeq} (b).

Conversely, let
\begin{equation*}
 M = \xmx_n \pmatrix{\omx_\nu & a \S_\nu^T \cr \S_\nu b^T & Z} \xmx_n,
\end{equation*}
where
$a, b, Z$
have the properties stated in the theorem. Then
$M \in M_n \cap S_n$
by Theorems \ref{tSblock} and \ref{tMblock}, and
\begin{align}
 M &= \rc 2 \pmatrix{a \S_\nu^T \jmx_\nu + \jmx_\nu \S_\nu b^T + \jmx_\nu Z \jmx_\nu & -a \S_\nu^T + \jmx_\nu \S_\nu b^T \jmx_\nu - \jmx_\nu Z \cr
         \jmx_\nu a \S_\nu^T \jmx_\nu - \S_\nu b^T - Z \jmx_\nu & - \jmx_\nu a \S_\nu^T - \S_\nu b^T \jmx_\nu + Z \cr}
\nonumber\end{align}
is of the form (\ref{ePstruct}), as can be checked by a straightforward calculation.
\phantom{.}\hfill$\qed$\nobreak\par\noindent
%\par\medskip\noindent
%The case of odd
%$n$
%can be treated very similarly, but there are some changes because in this case
%$\S_n$
%is even (instead of odd) under the operation of
%$\jmx$
%and does not add up to 0.
%\begin{Theorem}
%Let
%$M$
%be a
%$2n \times 2n$
%matrix, with
%$n$
%odd. Then
%$M$
%is a most perfect square of weight
%$w$
%if and only if there are vectors
%$a,b \in \{1_n\}^\bot$
%with
%$\jmx a = a, \jmx b = b$
%and an
%$n \times n$
%weightless generally associated matrix
%$Z$
%satisfying
%$\xi^T Z \eta = 0$
%$(\xi, \eta \in \{\S_n\}^T)$
%such that
%\begin{align}
% M = X \pmatrix{2 w E & a \S_n^T \cr \S_n b^T & Z \cr} X.
%\nonumber\end{align}
%\end{Theorem}
%\par\medskip\noindent
%The proof is analogous to the previous, but now
%$\jmx \S_n = \S_n;$
%the line (\ref{exalt}) now becomes
%\begin{align}
% X \S_{2n} &= \rc{\sqrt 2} \pmatrix{I & \jmx \cr \jmx & -I} \pmatrix{\S_n \cr -\S_n}
% = \rc{\sqrt 2} \pmatrix{\S_n - \jmx \S_n \cr \jmx \S_n + \S_n}
% = \pmatrix{0 \cr \sqrt{2} \S_n};
%\nonumber\end{align}
%instead of (\ref{eabpar}) we now have
%\begin{align}
% \S_n a^T = \jmx \S_n a^T \jmx &= \S_n a^T \jmx.
%\nonumber\end{align}
%This gives rise to the different reflexion parity of
%$a,b$
%in Theorem \ref{tmpblodd}.
%Note that here the off-diagonal blocks
%$V^T, W$
%will not be semi-magic (unless they are the null matrix).
\par\medskip\noindent
%[{\it Remark. \/}
%It would be interesting to see how the block representation elements
%$a, b, Z$
%of Theorems \ref{tMPSblock} and \ref{tmpblodd} are related to the block representation elements
%$a, b$
%of the unique corresponding reversible square, see Theorem \ref{tSally}.]
\par\medskip\noindent
The block representation of Theorem \ref{tMPSblock} can be used as a simple method of
constructing most perfect squares, as illustrated in the example
\begin{align}
 2 \xmx_6 \pmatrix{0  &  0  &  0  &  1  & -1  &  1 \cr
           0  &  0  &  0  & -2  &  2  & -2 \cr
           0  &  0  &  0  &  1  & -1  &  1 \cr
          -2  &  4  & -2  &  1  &  0  & -1 \cr
           2  & -4  &  2  & -1  &  0  &  1 \cr
          -2  &  4  & -2  &  1  &  0  & -1 \cr} \xmx_6 &=
 \pmatrix{-2  &  3  &  0  & -4  &  5  & -2 \cr
           1  & -2  & -1  &  5  & -6  &  3 \cr
          -2  &  3  &  0  & -4  &  5  & -2 \cr
           4  & -5  &  2  &  2  & -3  &  0 \cr
          -5  &  6  & -3  & -1  &  2  &  1 \cr
           4  & -5  &  2  &  2  & -3  &  0 \cr}.
\nonumber\end{align}
However, it turns out that the structure and construction of most perfect square matrices
is even more simple.
Indeed, they can be conveniently characterised, even without the use of the block
representation, in the following way.

\begin{Theorem}
\label{tMPSstruct}
A matrix
$M\in\mathbb{R}^{n \times n}$, $n = 2 \nu$ even,
is an element of $MPS_n$ if and only if
\begin{equation}
 M = \gamma \S_{n}^T + \S_{n} \delta^T,
\label{empstruct}\end{equation}
where
\begin{description}
\item{(a)}
in case
$\nu$
is even,
${\displaystyle\quad
 \gamma = \pmatrix{\tilde \gamma \cr - \tilde \gamma},
\quad
 \delta = \pmatrix{\tilde \delta \cr - \tilde \delta}, \quad \tilde\gamma, \tilde \delta \in {\mathbb R}^\nu,
}$
\item{(b)}
in case
$\nu$
is odd,
${\displaystyle\quad
 \gamma = \pmatrix{\tilde \gamma \cr \tilde \gamma},
\quad
 \delta = \pmatrix{\tilde \delta \cr \tilde \delta}, \quad \tilde\gamma, \tilde \delta \in \{1_\nu\}^T \subset {\mathbb R}^\nu.
}$
\end{description}
\par\medskip\noindent
The vectors
$\gamma, \delta$
can be obtained from
$M$
as
 $\gamma = \rc{n} M \S_{n}$, $\delta = \rc{n} M^T \S_{n}.$
\end{Theorem}
\par\medskip\noindent
{\em Proof.}
If $M \in MPS_n$, then it follows from the second equation in Theorem \ref{tsymeq}
(f) that
$\gamma := \frac{1}{n} M \S_n \in \{\S_n\}^\bot$.
\par
Let $v \in \{\S_n\}^\bot$; then, by the first equation in Theorem \ref{tsymeq} (f),
$M v \in \{\S_n\}^{\bot\bot} = \mathbb{R}\S_n$. Hence
$M v = f(v)\S_n$ $(v \in \{\S_n\}^\bot)$, where $f$ is a linear form on
$\{\S_n\}^\bot$.
By the Riesz representation theorem, there is a vector $\delta\in\{\S_n\}^\bot$
such that $f(v) = \delta^T v$, so $M v = \S_n \delta^T v$ $(v \in \{\S_n\}^\bot)$.
\par
Now any
$x \in \mathbb{R}^n$ can be written in the form
$x = \alpha \S_n + v$, with $\alpha\in\mathbb{R}$ and $v \in \{\S_n\}^\bot$; then
\begin{align*}
 M x &= \alpha M \S_n + M v
 = \alpha n \gamma + \S_n \delta^T v
 = (\gamma \S_n^T + \S_n \delta^T)(\alpha \S_n + v)
\\
 &= (\gamma \S_n^T + \S_n \delta^T) x,
\end{align*}
showing that $M$ is of the form (\ref{empstruct}).
\par
Writing
$\gamma = \pmatrix{\gamma_1 \cr \gamma_2}$,
$\delta = \pmatrix{\delta_1 \cr \delta_2}$,
with $\gamma_1, \gamma_2, \delta_1, \delta_2 \in \mathbb{R}^\nu$,
we find
\begin{equation*}
M = \pmatrix{\gamma_1\cr\gamma_2}\pmatrix{\S_\nu\cr\pm\S_\nu}^T
   + \pmatrix{\S_\nu\cr\pm\S_\nu}\pmatrix{\delta_1\cr\delta_2}^T
 = \pmatrix{\gamma_1\S_\nu^T + \S_\nu \delta_1^T & \pm \gamma_1\S_\nu^T + \S_\nu \delta_2^T \cr
   \gamma_2\S_\nu^T \pm \S_\nu \delta_1^T & \pm \gamma_2\S_\nu^T \pm \S_\nu \delta_2^T }.
\end{equation*}
Viewed in conjunction with (\ref{ePstruct}), this implies
$(\gamma_1 \pm \gamma_2)\S_\nu^T + \S_\nu (\delta_1 \pm \delta_2)^T = \omx_\nu$.
As
\begin{equation*}
 0 = \S_n^T \gamma = \pmatrix{\S_\nu \cr \pm\S_\nu}^T \pmatrix{\gamma_1 \cr \gamma_2}
 = \S_\nu^T (\gamma_1 \pm \gamma_2),
\end{equation*}
it follows that
\begin{equation*}
 0_\nu = \S_\nu^T (\gamma_1 \pm \gamma_2)\S_\nu^T + \S_\nu^T \S_\nu (\delta_1 \pm \delta_2)^T
 = 0_\nu + \nu (\delta_1 \pm \delta_2)^T,
\end{equation*}
so $\delta_2 = \mp \delta_1$. An analogous calculation gives $\gamma_2 = \mp \gamma_1$.
\par
By Lemma \ref{lMNSVeigen} and Theorem \ref{tsymeq} (a), (f) (with $u = v = 1_n$),
$1_n$ is an eigenvector with eigenvalue $0$ for both
$M$ and $M^T$. Since $\S_n^T 1_n = 0$, it follows that $\gamma, \delta \in \{1_n\}^\bot$.
If $\nu$ is even, this will be satisfied for any $\gamma_1$, $\delta_1$ in view of
the above structure; if $\nu$ is odd, it gives the further condition that
$\gamma_1, \delta_1$ are orthogonal to $1_\nu$.
\par
For the converse, a straightforward calculation shows that any $M$ of the form
(\ref{empstruct}), with $\gamma$, $\delta$ satisfying the hypotheses, has the properties
listed in Theorem \ref{tsymeq} (a), (f) and (\ref{ePstruct}).
\phantom{.}\hfill$\qed$\nobreak\par\noindent
\par
\bigskip\noindent
The two terms in the representation (\ref{empstruct}) are obviously rank 1
matrices (if non-null), so we can immediately draw the following conclusion.
\begin{Corollary}
A most perfect square matrix has at most rank 2 if its weight is 0, at most rank 3 in general.
\end{Corollary}
\par\medskip\noindent
As a further consequence, we find an equivalent criterion for {\em parasymmetry}
of a most perfect (and in particular magic) square matrix $M$, defined in \cite{hill}
as symmetry of its square $M^2$.
Indeed, as $\S_n^T \gamma = \delta^T \S_n = 0$ and therefore
$M^2 = 2n \gamma \delta^T + (\delta^T \gamma) \S_n \S_n^T$, the following is evident.
\begin{Corollary}
A weight $0$ most perfect square
\begin{equation*}
 M = \gamma \S_n^T + \S_n \delta^T \in MPS_n
\end{equation*}
is parasymmetric if and only if $\gamma, \delta$ are linearly dependent.
\end{Corollary}

\medskip\noindent
By Lemma \ref{lalginters}, $MPS_n$ is the complement to a product symmetry
type $NQS_n := N_n \cap Q_n \cap S_n$, so that $NQS_n \oplus MPS_n$ is again
a $\mathbb{Z}_2$-graded algebra.

NQS-type matrices have the following block representation.

\begin{Theorem}\label{tNQSblock}
Let
$M\in\mathbb{R}^{n \times n}$, $n = 2 \nu$ even.
Then $M\in NQS_n$
if and only if
\begin{align}
 M &= \xmx_n \pmatrix{Y & V^T \cr W & Z} \xmx_n,
\label{eNQSblock}\end{align}
where
$Y\in B_\nu \cap S_\nu$,
$Z\in B_\nu \cap N_\nu$,
and
$V, W\in A_\nu$
have the properties
\begin{equation}
\label{eVWcond}
V 1_\nu = W 1_\nu = 0_\nu, \qquad V^T \S_\nu = W^T \S_\nu = 0_\nu.
\end{equation}
\end{Theorem}

\noindent
The conditions on $V$ and $W$ in Theorem \ref{tNQSblock} mean that these matrices
have symmetry (A) with weight 0, and all their row sums and all alternating column
sums vanish.
\par\medskip\noindent
{\it Proof. \/}
If $M \in NQS_n$, then by combining the block representations of Theorems
\ref{tSblock} and \ref{tNblock}, we find that
(\ref{eNQSblock}) holds with $Y\in S_n$, $Z\in N_n$
and with
$V, W$
satisfying (\ref{eVWcond}).
As
$M$
also has property (Q), we find from (\ref{eQstruct}) that
\begin{align}
 \xmx_n M \xmx_n &= \pmatrix{A + B \jmx_\nu + \jmx_\nu B + \jmx_\nu A \jmx_\nu & A \jmx_\nu - B + \jmx_\nu B \jmx_\nu - \jmx_\nu A \cr
                  \jmx_\nu A + \jmx_\nu B \jmx_\nu - B - A \jmx_\nu & \jmx_\nu A \jmx_\nu - \jmx_\nu B - B \jmx_\nu + A}
\nonumber\end{align}
and can read off, upon multiplication with $\jmx_\nu$ on both sides, that
$V^T, W \in A_\nu$
and
$Y, Z\in B_\nu$.
\par
Conversely, assume
$M$
is given by (\ref{eNQSblock}), where
$V,W,Y,Z$
have the stated properties. Then
$M$
is semimagic by Theorems \ref{tSblock} and \ref{tNblock}, and
\begin{align}
 M &= \rc 2 \pmatrix{Y + V^T \jmx_\nu + \jmx_\nu W + \jmx_\nu Z \jmx_\nu & Y \jmx_\nu - V^T + \jmx_\nu W \jmx_\nu - \jmx_\nu Z \cr
                   \jmx_\nu Y + \jmx_\nu V^T \jmx_\nu - W - Z \jmx_\nu & \jmx_\nu Y \jmx_\nu - \jmx_\nu V^T - W \jmx_\nu + Z}
\nonumber\end{align}
has property (Q).
\phantom{.}\hfill$\qed$\par\noindent

\medskip\noindent
{\it Remark.\/}
The $\mathbb{Z}_2$-graded algebra $NQS_n \oplus MPS_n$ is a subalgebra of $S_n$,
but not all of $S_n$, since the types N and M are here linked to types Q and S,
respectively.

As the `odd' part of a $\mathbb{Z}_2$-graded algebra, $MPS_n$ is not itself a
subalgebra of the full matrix algebra; however, it has the property that the
product of any {\it three\/} elements of $MPS_n$ is again an element of $MPS_n$.
If, taking the formula (\ref{empstruct}) as a motivation, we introduce the notation
$(\gamma; \delta) := \gamma \S_n^T + \S_n \delta^T$
for most perfect square matrices,
then the triple product can be expressed as
\begin{align}
 (\gamma_1; \delta_1) (\gamma_2; \delta_2) (\gamma_3; \delta_3)
 &= n ((\delta_2^T \gamma_3) \gamma_1; (\delta_1^T \gamma_2) \delta_3).
\nonumber\end{align}

\section{Composite Symmetry: Reversible Squares}
\label{sCSRS}

\noindent
We now turn to reversible square matrices, defined as those which have symmetry
properties (R) and (V). Although the definition of these properties does not
refer to a weight $w$, it turns out that reversible squares always have the
associated symmetry property (A) and hence a hidden weight.

\begin{Lemma}\label{lrevisga}
Any reversible square matrix has property (A) with some weight $w\in\mathbb{R}$.
\end{Lemma}
\par\medskip\noindent
{\it Proof. \/}
Let $M\in\mathbb{R}^{n\times n}$ be a reversible matrix, so $M\in R_n$ and
$M$ has property (V), which, following the proof of Theorem \ref{tsymeq} (e),
can be seen to be equivalent to
\begin{equation}
\label{eVpropeq}
u^T M v^T = 0 \qquad (u, v \in\{1_n\}^\bot).
\end{equation}
By Theorem \ref{tsymeq} (b), we only need to show that there is $w\in\mathbb{R}$
such that $M + \jmx_n M \jmx_n = 2 w \emx_n$.

Consider the two orthogonal projectors (symmetric idempotent matrices)
 $P = \rc 2 (\imx_n + \jmx_n)$, $Q = \rc 2 (\imx_n - \jmx_n)$;
clearly
$P^T = P, Q^T = Q$
and
$P + Q = \imx_n.$
Moreover,
$\jmx_n P = P = P \jmx_n$
and
$\jmx_n Q = -Q = Q \jmx_n.$
Also,
$P 1_{n} = 1_{n}$
and
$Q 1_{n} = 0_n.$
Using these properties, we deduce
\begin{align}
 M &+ \jmx_n M \jmx_n = (P + Q)^T (M + \jmx_n M \jmx_n) (P + Q)
\nonumber\\
 &= P M P + P \jmx_n M \jmx_n P + P M Q + P \jmx_n M \jmx_n Q + Q M P + Q \jmx_n M \jmx_n P + Q M Q + Q \jmx_n M \jmx_n Q
\nonumber\\
 &= 2 P M P + 2 Q M Q
 = 2 P M P,
\nonumber\end{align}
observing in the last step that, by (\ref{eVpropeq}),
$u Q M Q v = 0$
for all
$u, v \in {\mathbb R}^{n},$
since
$1_{n}^T Q u = (Q 1_{n})^T u = 0$
and similarly for
$v.$
By analogous reasoning, if
$u \in \{1_{n}\}^\bot,$
then also
$P u \in \{1_{n}\}^\bot$,
so
\begin{align}
 (M + \jmx_n M \jmx_n)u = 2 P M P u &= P M P u + P M (\jmx_n P) u = P (M + M \jmx_n) P u
 = \omx_n
\nonumber\end{align}
by Theorem \ref{tsymeq} (d).
Thus the dimension of the kernel of
$M + \jmx M \jmx$
is at least
$n-1,$
so this matrix has at most rank 1.
If it has rank 0, then $M \in A_n$ by Theorem \ref{tsymeq} (b).

Assuming rank 1 in the following,
we can rewrite $M + \jmx_n M \jmx_n = P (M + \jmx_n M) P$ as above;
by Theorem \ref{tsymeq} (d), the range of
$M + \jmx_n M$
is $\mathbb{R} 1_n$,
which is invariant under the action of
$P.$
Thus, the range of the rank 1 matrix
$M + \jmx_n M \jmx_n$
is $\mathbb{R} 1_n$.
\par
In summary, both the kernel and the range of the rank 1 matrix
$M + \jmx_n M \jmx_n$
are equal to the kernel and range of $\emx_n$, respectively.
As a rank 1 matrix is determined up to a multiplicative constant by its
kernel and range, it follows that
$M + \jmx_n M \jmx_n = 2 w \emx_n$
for some
$w\in\mathbb{R}.$
\phantom{.}\hfill$\qed$\par\noindent

\bigskip
\noindent
The set $V_n$ has, in addition to (V), the defining requirement that the sum of
all matrix entries vanish. For a matrix with property (A) this is the case if
and only if the weight $w$ vanishes. Hence we can see that any reversible square
matrix is a sum of an element of $R_n \cap V_n$ and a multiple of $\emx_n$
(this decomposition being unique since $\emx_n \notin V_n$).

Therefore it makes sense to focus on the space of {\it weightless reversible
squares\/} $RV_n := R_n \cap V_n$. Its elements have the following block representation.

\begin{Theorem}
\label{tRVblock}
Let
$M \in \mathbb{R}^{n\times n}$.
Then
$M$
is an element of $RV_n$
if and only if there are vectors
$a, b \in {\mathbb R}^\nu$
such that
\begin{align}
 M &= \xmx_n \pmatrix{\omx_\nu & 1_\nu a^T \cr b 1_\nu^T & \omx_\nu} \xmx_n
\label{eSally}\end{align}
if $n = 2\nu$ is even,
\begin{equation*}
 M = \xmx_{n} \pmatrix{\omx_\nu & 0_\nu & \sqrt 2\,1_\nu a^T \cr 0_\nu^T & 0 &  a^T
     \cr \sqrt 2\,b 1_\nu^T &  b & \omx_\nu} \xmx_{n}
\nonumber\end{equation*}
if $n = 2\nu + 1$ is odd.
\end{Theorem}
\par\medskip\noindent
{\it Proof.\/}
In the even-dimensional case, the result follows from comparison of the block representations in
Theorems \ref{tRblock} and \ref{tVblock}, noting that $\emx_\nu\notin V_\nu$
and therefore $\gamma = 0$ in (\ref{eRbleven}).
Similarly, in the odd-dimensional case,
comparison of (\ref{eRblodd}) with (\ref{eVblodd}) shows that
$v = y = \gamma 1_\nu$, and hence, in the central matrix entry,
$\frac{4\nu}{2\nu-1}\,\gamma = \gamma$,
which forces $\gamma = 0$ (as does the identity of the top left $\nu\times\nu$
blocks).
\phantom{.}\hfill$\qed$

\bigskip
\noindent
{\it Remarks.\/}
1.
From the formulae in Theorem \ref{tRVblock}, the block representation of a general
reversible square can be obtained by adding a suitable multiple of the block
representation of $\emx_n$; this gives
\begin{equation*}
 M = \xmx_n \pmatrix{2 w \emx_\nu & 1_\nu a^T \cr b 1_\nu^T & \omx_\nu} \xmx_n
\end{equation*}
if $n = 2\nu$ is even,
\begin{equation*}
 M = \xmx_{n} \pmatrix{2 w \emx_\nu & \sqrt 2 w\,1_\nu &\sqrt 2\,1_\nu a^T \cr \sqrt 2 w\,1_\nu^T & w &  a^T
     \cr \sqrt 2\,b 1_\nu^T &  b & \omx_\nu} \xmx_{n}
\nonumber\end{equation*}
if $n = 2\nu + 1$ is odd.

2.
Comparison of the block structures in Theorem \ref{tRVblock} and Lemma \ref{ttopsy}
shows that type RV matrices are type A matrices with rank 1 blocks with
constant rows and columns, respectively.
It follows that reversible square matrices always have rank 2 or lower.

\bigskip\noindent
The block representation of Theorem \ref{tRVblock} reveals a further connection
of the space $RV_n$ with the space $AS_n := A_n \cap S_n$ of (weightless)
associated magic square matrices.
\begin{Corollary}
\label{cRVisAV}
Let $n \in \mathbb{N}$.
Then $RV_n = AV_n$. Consequently, $A_n = RV_n \oplus AS_n$.
\end{Corollary}
\par\medskip\noindent
{\it Proof. \/}
From the block representations in Lemma \ref{ttopsy} and Theorem \ref{tVblock}, we
see that matrices in $AV_n$ have the same block representation as those in $RV_n$.
The second statement follows from $S_n \oplus V_n = \mathbb{R}^{n\times n}$.
\phantom{.}\hfill$\qed$

\bigskip\noindent
Corollary \ref{cRVisAV} shows that reversible square matrices can equivalently be defined
as square matrices with the properties (A) and (V); as $\emx_n$ has both of these
properties, this includes general weighted reversible square matrices.

The space $RV_n$ does not form a matrix algebra by itself; however, as we know
that $S_n \oplus V_n$ is a $\mathbb{Z}_2$-graded algebra by Theorem \ref{tSValg},
and $R_n$ is an algebra by Theorem \ref{tRalg}, it follows by Lemma \ref{lalginters}
(a) that $R_n = RS_n \oplus RV_n$ (where $RS_n := R_n \cap S_n$) is a
$\mathbb{Z}_2$-graded algebra.

\bigskip\noindent
{\it Remark.\/}
From Theorems \ref{tSblock} and \ref{tRblock} it is apparent that the elements
of $RS_n$ have the block representation
\begin{equation*}
 M = \xmx_n \pmatrix{\gamma\emx_\nu & \omx_\nu \cr \omx_\nu & Z} \xmx_n
\end{equation*}
in the even-dimensional case and
\begin{equation*}
 M = \xmx_n \pmatrix{2\gamma\emx_\nu & \sqrt 2\,\gamma 1_\nu & \omx_\nu \cr
 \sqrt 2\,\gamma 1_\nu^T & \gamma & 0_\nu^T \cr \omx_\nu & 0_\nu & Z} \xmx_n
\end{equation*}
in the odd-dimensional case; equivalently, they are the sum of a multiple of
$\emx_n$ and of a matrix of form
\begin{equation*}
 \xmx_n \pmatrix{\omx_{n-\nu} & \omx \cr \omx & Z} \xmx_n \in B_n \cap S_n
\end{equation*}
with arbitrary $Z \in \mathbb{R}^{\nu\times\nu}$.

While this shows that $RS_n$ is only a proper subspace of $BS_n := B_n \cap S_n$,
we note, applying Lemma \ref{lalginters} (b) to $B_n \oplus A_n$ and
$S_n \oplus V_n$, that $BS_n \oplus RV_n$ is another $\mathbb{Z}_2$-graded
algebra; in particular, the product of a weightless reversible square and a
balanced semimagic square matrix is a weightless reversible square matrix.
The latter statement clearly extends to general, weighted reversible square
matrices.

\bigskip\noindent
{\bf Acknowledgement.\/}
Sally Hill's research was funded by a Cardiff University EPSRC DTP grant EP/L504749/1.

\small{School of Mathematics\\
Cardiff University\\
Senghennydd Road\\
Cardiff CF24 4AG\\
UK\\
email: LettingtonMC@cardiff.ac.uk}; SchmidtKM@cardiff.ac.uk
\bigskip
\end{document}